\documentclass[preprint]{elsarticle}
\usepackage{graphicx} 
\usepackage{amsmath}
\usepackage{amsthm}
\usepackage{amssymb}
\usepackage{faktor}
\usepackage{mathtools}
\usepackage{tikz-cd}
\usepackage{subcaption}
\usepackage{url}
\usepackage{float}

\usepackage{qting}

\newtheorem{remark}{Remark}
\newtheorem{proposition}{Proposition}
\newtheorem{lemma}{Lemma}
\newcommand{\bigpartialderiv}[2]{ \frac{\partial {#1}}{\partial {#2} } }

\cortext[cor]{Corresponding author}
\affiliation[1]{organization={University of Cologne}, country={Germany}}
\affiliation[2]{organization={Universidad Politécnica de Madrid}, country={Spain}}
\affiliation[3]{organization={Max Planck Institute for Plasma Physics}, country={Germany}}
\affiliation[4]{organization={Technical University of Munich}, country={Germany}}
\affiliation[5]{organization={University of Waterloo}, country={Canada}}

\begin{document}

\title{SBP-FDEC: Summation-by-Parts Finite Difference Exterior Calculus}

\author[1]{Daniel Bach\corref{cor}}
\ead{daniel.bach@uni-koeln.de}
\author[2]{Andrés M. Rueda-Ramírez}
\ead{am.rueda@upm.es}
\author[3,4]{Eric Sonnendrücker}
\ead{eric.Sonnendruecker@ipp.mpg.de}
\author[5]{David C. Del Rey Fernández}
\ead{ddelreyfernandez@uwaterloo.ca}
\author[1]{Gregor J. Gassner}
\ead{ggassner@uni-koeln.de}

\begin{abstract}
We demonstrate that we can carry over the strategy of Finite Element Exterior Calculus (FEEC) to Summation-by-Parts (SBP) Finite Difference (FD) methods to achieve divergence- and curl-free discretizations. This is not obvious at first sight, as for SBP-FD no basis functions are known, but only values and derivatives at points. The key is a remarkable analytic relationship that enables us to construct compatible operators using integral and nodal degrees of freedom. Pre-existing SBP-FD matrix operators can then be used to obtain nodal values from the integral degrees of freedom to derive a scheme with the desired properties.   
\end{abstract}


\maketitle

\section{Introduction}

The theory of Finite Element Exterior Calculus (FEEC) by Arnold, Falk, and Winther \cite{arnold2006,arnold2010}, Arnold \cite{arnold2018}, with foundational contributions by Monk \cite{monk_2003}, is a state-of-the-art finite element framework that uses the structure of the de Rham cohomology to preserve topological structures, such as the cohomology groups of the region, or manifold, on which the Partial Differential Equations (PDEs) in question are posed. 
The primary example used in these foundational papers is the Hodge-Laplacian, but another major application of this framework are problems where an exactly divergence- or curl-free solution is needed.
For instance, FEEC has been used to satisfy the divergence-free condition of the magnetic field when discretizing the time-dependent Maxwell's equations (see, for instance, \cite{arnold2018} and \cite{monk_2003}) and Vlasov--Maxwell equations (see, for instance \cite{kraus_2017,Campos-Pinto2022Variational}).
The theory has also been applied to other PDEs that need to satisfy divergence constraints, such as the incompressible Navier--Stokes equations.
See, for example, Palha et al. \cite{palha_2017} or Carlier et al. \cite{carlier_2024}.
In FEEC, compatible ansatz spaces and projections are used to discretize the PDEs.
For that purpose, two (in 1D), three (in 2D), or four (in 3D) discrete spaces and conforming projections onto these spaces are necessary to discretize the de Rham cochain complex.
The FEEC framework has also recently been extended to broken spaces, such that only local mass matrices are needed \cite{gucclu2023broken}.

In the literature, a variety of different ansatz spaces and projections have been proposed to discretize the de Rham cochain complex, most of which are polynomial or piecewise polynomial in nature.
Classical examples include Nédélec spaces of the first and second kind \cite{nedelec1980,nedelec1986}, Whitney forms \cite{whitney1957} or Raviart--Thomas elements \cite{raviart_thomas_1977}.
These predate the formulation of the FEEC framework.
More recently, spline elements using a tensor-product ansatz on hexahedral spaces have also been used, for example, in Buffa et al. \cite{buffa_2011} and Evans et al. \cite{evans_2013}.

These approaches differ not only in the choice of finite element spaces, but also in the projection operators they employ, which leads to variations in the associated degrees of freedoms and basis functions of each approach.
For instance, the approach of Gerritsma \cite{gerritsma2011}  and Gerritsma et al. \cite{preprint_gerritsma2011}, which motivated the present work, uses nodal and integral degrees of freedom on a one-dimensional subgrid per element.
Using tensor products, this construction extends to nodal, line-integral, area-integral, and volume-integral degrees of freedom on a multi-dimensional level.
Other approaches, as in Nédélec \cite{nedelec1980,nedelec1986}, use higher moment degrees of freedom instead of a subgrid approach.

There is a rich literature of discrete operators on finite-difference-like degrees of freedom using discrete derivative operators to mimic the de Rham complex discretely to achieve similar structure preserving properties as FEEC.
These schemes are also known as mimetic finite-difference methods.
One survey paper by Lipnikov et al. \cite{lipnikov_2014} gives a comprehensive overview of these methods.
An early example of such schemes is given in the book of Samarskii \cite{samarskii_translation}.
Another important work is the method of Yee \cite{yee_1966} for Maxwell's equations, which is the foundation for the class of finite-difference time-domain methods. 
It uses a staggered space-time grid for the different fields, which automatically preserves the divergence condition of the magnetic field.
Another approach of Hyman and Shashkov in \cite{hyman_1997_2,hyman_1997} and \cite{hyman_1998} uses the derivation of discrete analogues for the \textbf{grad}, \textbf{curl}, and div operators and their adjoints.
In addition to the natural operators, which are these discrete analogues,  arising in this context, mimetic discretizations rely on a reconstruction operator to compute point values of the different fields involved. For high-order mimetic finite differences, these are mostly interpolation operators based on the degrees of freedoms in the neighboring cells, as for example in high-order finite volume reconstructions. In \cite{Kormann2024A-Dual-Grid} a mimetic method on staggered grids is developed for the Vlasov--Maxwell equations, where the reconstructions are defined by a sliding stencil of Lagrange interpolation or histopolation using the neighboring cells. 

A different class of finite-difference methods, called Summation-by-Parts Finite Differences (SBP-FD), introduced Kreiss and Scherer \cite{kreiss_1974} with key theoretical contributions from Strand \cite{strand1994}, were developed to re-create energy estimates of classic Finite-Element Methods (FEM).
SBP operators use a discrete norm (mass) matrix together with a derivative matrix that operate on nodal degrees of freedom.
These operators are constructed such that they mimic integration by parts discretely, enabling the derivation of discrete energy estimates, see Olsson in \cite{olsson_sbp_1, olsson_sbp_2}, and entropy estimates, see for instance Crean et al. \cite{crean_2018} or Carpenter et al. \cite{carpenter_2016}.
A defining feature of these methods is the equivalence of their strong and weak formulations, which is central to obtaining stability results.

These schemes can be extended to a multi-element formulation with the Simultaneous Approxmation Term (SAT) approach, treating interface and boundary conditions weakly. 
For further information, see, for instance, the review papers \cite{svard_2014, fernandez2014review}.
SBP operators for the first and second derivative have also been developed for non-polynomial approximation spaces, for instance for trigonometric functions or radial basis functions, see Glaubitz et al. \cite{glaubnitz_2023, glaubnitz_2024}.
The SBP framework has also influenced methods beyond finite-difference schemes.
One example of this is the Discontinuous Galerkin Spectral Element Method (DGSEM) using a Legendre--Gauss--Lobatto (LGL) grid, see Black \cite{black1999conservative} and Kopriva \cite{kopriva2006}, where the mass and derivative matrices form an SBP operator, see Gassner \cite{gassner2014}.

Because classical SBP-FD methods rely solely on nodal degrees of freedom and a discrete first-order derivative operator acting on those nodes, it is not immediately clear how to extend them analogously to FEEC to achieve similar structure-preserving properties. 
Motivated by the analogy between classic nodal FEM and SBP-FD schemes used to mimic the energy analysis, the central question of this paper is whether it is possible to develop a Summation-by-Parts Finite Difference Exterior Calculus (SBP-FDEC) analogous to FEEC. The major scientific obstacle to overcome lies in the fact that, when using only nodal DOFs (SBP-FD), without known analytic basis functions, it is not obvious how to derive discrete derivative operators and projections that satisfy a discrete de Rham cohomology. 

To address this issue, we exploit important connections between DGSEM-LGL and SBP-FD, and borrow important ingredients from the work of Gerritsma \cite{preprint_gerritsma2011} and Gerritsma et al. \cite{gerritsma2011}, who developed FEEC operators using LGL spectral element operators. In Gerritsma et al.'s approach, a subgrid structure is imposed on each LGL element, accompanied by two compatible polynomial spaces and two projection operators in 1D: (i) the polynomial interpolation space and (ii) the polynomial \textit{histopolation} space.  

In this paper, we show that, although the histopolation basis polynomials are derived from the LGL Lagrange basis functions, their nodal values depend solely on the (discrete) nodal derivative values of the Lagrange basis.
This important observation and the close relationship between LGL spectral element operators and SBP-FD leads us to the central idea of this work: it is possible to use nodal values like classical SBP-FD operators, while simultaneously accessing sub-grid interval DOFs like in the work of Gerritsma et al.
The interval DOFs can be converted back into nodal values by using the derivative matrix of the standard SBP-FD derivative operator, which can be accomplished without requiring explicit basis functions.

This strategy enables us to invoke the fundamental theorem of calculus to construct a discrete complex that is compatible with the de Rham complex, while still performing numerical computations with existing SBP-FD operators using nodal values derived from the integral DOFs.
The extension to 2D and 3D via a tensor product ansatz is then straightforward, and when the SBP-FD operator grid contains the boundary nodes, it is also easy to extend the method to multiple elements.
Overall, we are able to derive the class of SBP-FDEC schemes.

The remainder of the paper is organized as follows. 
In section 2, we revisit and generalize the ansatz of Gerritsma \cite{gerritsma2011} to the SBP-FD framework.
We include extensions to multiple elements and derivations of discrete 2D and 3D derivative operators based on existing SBP-FD operators.
We further show that the cohomology spaces of the de Rham complexes for non-periodic rectangular or hexahedral domains without imposed boundary conditions are preserved by the discrete complexes.
In section 3, we apply this theory to the 2D transverse electric homogeneous Maxwell's equation to demonstrate the capability of the SBP-FDEC to  exactly preserve the divergence-free electric field.
In section 4, we derive the weak and strong form of the SBP-FDEC method for the Maxwell's equations on a periodic domain and prove the discrete divergence-free property and semi-discrete energy conservation, which can then be extended to fully discrete energy conservation using appropriate established time-integration methods. 
In section 5, we numerically verify our method both on a periodic and a non-periodic test case, showing the experimental order of convergence fitting to the established SBP-FD operators, demonstrating that the expected convergence rate is obtained as well as that the scheme is discretely divergence-free and energy preserving.
The last section concludes the paper.

\section{Mimetic Finite Difference Operators}

\subsection{1D Basis Functions}

\label{subsection:basis_functions}

Our goal in this section is to construct discrete spaces $U_l$ and $U_h$ and compatible projections $p^0$ and $p^1$ to obtain a discrete 1D de Rham complex, so that the diagram of Fig. \ref{fig:1d_discrete_de_rham} commutes.
This means, that we get the same result if we first apply the derivative operator $\bigpartialderiv{}{x}$ and then $p^1$ or first $p^0$ and then $\bigpartialderiv{}{x}$, and we also get the same result if we first apply $\bigpartialderiv{}{x}$ and then $R^1$ or first $R^0$ and then $\underline{\Delta}$, where $R^0$ and $R^1$ are reduction operators, mapping functions to the vector of their DOFs.
The discrete de Rham complex {ensures} that our spaces are compatible with the derivative operator $\bigpartialderiv{}{x}$.
This we can later extend with a tensor product ansatz to 2D and 3D to obtain spaces and projections compatible with the $\textbf{curl}$ and $\text{div}$ operators, to achieve divergence-free discretizations.

To obtain such a discrete complex, we will generalize the ansatz proposed by Gerritsma \cite{gerritsma2011}, who derives two compatible polynomial spaces using nodal and integral degrees of freedom on a given grid.
Instead of using polynomial spaces, however, we assume an unknown general Lagrange basis that we only know through discrete nodal function values and the corresponding discrete nodal derivative values. This assumption allows us to connect nodal function values to a new set of discrete integral degrees of freedom and further derive corresponding discrete representations and matrix operators. 

With this, we are able to derive a discrete de Rham complex, shown in Fig. \ref{fig:1d_discrete_de_rham}, on the spaces of the nodal and integral DOFs respectively.
\begin{figure}[h]
    \centering
    \begin{tikzpicture}
        \label{derham1d_feec}
       	\node (1) at (2,0) {$H^1$};
       	\node (2) at (4,0) {$L^2$};
       	\node (3) at (2,-1.5) {$U_l$};
       	\node (4) at (4,-1.5) {$U_h$};
       	\node (5) at (2,-3) {$\mathbb{R}^{N+1}$};
       	\node (6) at (4,-3) {$\mathbb{R}^{N}$};
       	\draw[->] (1) -- node [midway, above] {$\bigpartialderiv{}{x}$} (2);
       	\draw[->] (3) -- node [midway, above] {$\bigpartialderiv{}{x}$} (4);
       	\draw[->] (5) -- node [midway, above] {$\underline{\Delta}$} (6);
       	\draw[->] (3) -- node [midway, right] {$R^0$} (5);
       	\draw[->] (4) -- node [midway, right] {$R^1$} (6);
       	\draw[->] (1) -- node [midway, right] {$p^0$} (3);
       	\draw[->] (2) -- node [midway, right] {$p^1$} (4);
    \end{tikzpicture}
    \caption{Discrete de Rham complex in 1D.}
    \label{fig:1d_discrete_de_rham}
\end{figure}

Let $\{x_i\}_{i = 0}^N$ be a grid in the interval $I$ with $x_0 = a$ and $x_N = b$ and let $\{l_i\}_{i = 0}^N$ be functions in $C^1\left([a,b]\right)$, which fulfill the Lagrange property for the grid points,
\begin{equation}
    \label{eq:Lagrange_property}
    l_i\left(x_j\right) = \delta_{i,j} \qquad \forall i,j \in \{0,\dots,N\},
\end{equation}
where $\delta_{i,j}$ is the Kronecker delta.
In addition, let $U_l$ be the span of the Lagrange functions $\{l_i\}_{i = 0}^N$. We demand that both the constant and linear functions are contained in $U_l$ to ensure that both $U_l$ and its derivative space contain the constant functions.
This is equivalent to the conditions
\begin{align}
    \label{eq:Lagrange_constant_functions}
    \sum_{i = 0}^N l_i(x) &= 1, \qquad \forall x \in [a,b], \\
    \label{eq:Lagrange_linear_functions}
    \sum_{i = 0}^N x_i l_i(x) &= x, \qquad \forall x \in [a,b].   
\end{align}
Since we need to calculate the derivatives, we want a basis representation for $U_h := \bigpartialderiv{}{x}\left(U_l\right)$.
For this, let $\{h_i\}_{i=1}^N$ be the functions in $C^0\left([a,b]\right)$ defined by
\begin{equation}
    \label{eq:histopolation_function}
    h_i = -\sum_{j = 0}^{i-1} \bigpartialderiv{l_j}{x}, \qquad \forall i \in \{1,\dots,N\}.
\end{equation}
We call these functions edge functions or histopolation functions associated with the Lagrange functions $\{l_i\}_{i = 0}^N$, using the terminology of Gerritsma \cite{gerritsma2011}.
It is straightforward to check that the histopolation property
\begin{equation}
    \label{eq:hist_int_def}
    \int_{x_{i-1}}^{x_{i}} h_j(x) \text{ d}x = \delta_{i,j}, \qquad \forall i,j \in \{1,\dots,N\},
\end{equation}
holds.
This property is analogous to the Lagrange property on the grid nodes but instead applies to the integrals over the sub-intervals defined by these nodes. 
It follows that $\{h_i\}_{i=1}^N$ are linearly independent functions and the span of these functions is a space of dimension $N$.
Since the $h_i$ are sums of derivatives of the Lagrange functions, we know that $h_i \in U_h$ for all $i \in \{1,\dots,N\}$.
On the other hand, the dimension of $U_l$ is $N+1$ and $U_l$ contains all constant functions.
Hence, $\text{dim}\left(U_h\right) = \text{dim}\left(\bigpartialderiv{}{x}\left(U_l\right)\right) = N$, which shows that $U_h$ is the span of the histopolation functions.
Given a function $g \in U_h$, the histopolation property \eqref{eq:hist_int_def} of the edge functions implies that the functions $h_i$ form a complete basis for $U_h$, i.e.
\begin{equation}
    \label{eq:histopolation}
    g(x) = \sum_{i = 1}^N\left(\int_{x_{i-1}}^{x_{i}} g(s) \text{ d}s\right) h_i(x) =: \sum_{i = 1}^N\overline{g}_i h_i(x).
\end{equation}
In addition, we know that all constant functions are contained in $U_h$ since the linear functions are contained in $U_l$. For $g(x) \equiv 1$ we have due to \eqref{eq:histopolation}
\begin{equation}
    g(x) = \sum_{i = 1}^N\left(\int_{x_{i-1}}^{x_{i}} 1 \text{ d}s\right) h_i(x) = \sum_{i = 1}^N \left(x_{i}-x_{i-1}\right) h_i(x) = 1, \qquad \forall x \in [a,b]. \\
\end{equation}
Let the interpolation operator $p^0$ and the histopolation operator $p^1$ be defined as
\begin{align}
    \label{eq:def_interpolation}
    p^0(f) &:= \sum_{i=0}^N f_i l_i = \sum_{i=0}^N f(x_i) l_i, \qquad &f \in C^1(I), \\
    \label{eq:def_histopolation}
    p^1(g) &:= \sum_{i = 1}^N \overline{g}_i h_i = \sum_{i = 1}^N\left(\int_{x_{i-1}}^{x_{i}} g(s) \text{ d}s\right) h_i, \qquad &g \in C^0(I),
\end{align}
We have for any function $f \in C^1(I)$ that
\begin{equation}
    \label{eq:fundamental_theorem}
    \int_{x_{i-1}}^{x_i} \bigpartialderiv{f}{x}(x)  \text{ d}x = f(x_i) - f(x_{i-1}), \qquad \forall i \in \{1,\dots,N\},
\end{equation}
and using Eq. \eqref{eq:histopolation_function}, we obtain the key relation
\begin{equation}
    \label{eq:derivative_formula}
p^1\left(\bigpartialderiv{}{x}f\right) = \sum_{j=1}^N \left(f_{j}-f_{j-1}\right)h_j = \bigpartialderiv{}{x}\left(\sum_{i=0}^N f_i l_i\right) = \bigpartialderiv{}{x}p^0(f),
\end{equation}
where we used
\begin{equation}
    -\sum_{i = 0}^{N-1}\bigpartialderiv{l_i}{x} = \bigpartialderiv{l_N}{x},
\end{equation}
which follows from Eq. \eqref{eq:Lagrange_constant_functions}.
If we have the derivative values of the $l_i$ functions on our grid points, we can define the derivative matrix $\underline{D} = (D_{i,j})_{i,j = 0}^N$ with entries
\begin{equation}
    \label{eq:derivative_matrix}
    D_{i,j} = \bigpartialderiv{l_j}{x}(x_i), \qquad \forall i,j \in \{0,\dots,N\},
\end{equation}
which acts as a discrete derivative operator
\begin{equation}
    \label{eq:discrete_derivative}
    \underline{f'} = \underline{D}\,\underline{f},
\end{equation}
where $\underline{f'}$ contains the nodal values of the approximate derivative of the function $f$. 
It is important to note that the nodal derivative matrix $\underline{D}$ allows us to directly compute the histopolation basis $h_i$ at the grid points using \eqref{eq:histopolation_function}, defining the histopolation Vandermonde matrix $\underline{V} \in \mathbb{R}^{N+1 \times N}$:
\begin{equation}
    \label{eq:vandermonde}
    V_{k,i} := h_i(x_k) = -\sum_{j=0}^{i-1}D_{k,j}, \qquad \forall i \in \{1,\dots,N\}, k \in \{0,\dots,N\}.
\end{equation}
Note that we number the row indices from $0$ to $N$ and the column indices from $1$ to $N$ to stay consistent with the numbering of the histopolation functions.
\begin{remark}
    For SBP-FD schemes, typically no continuous basis functions $l_i$ and/or $h_i$ are known. Only the value of the functions and their derivatives are known at the grid nodes. 
    Hence, the relation Eq. \eqref{eq:vandermonde} is a key observation that we use to construct discrete mimetic FD operators.
\end{remark}
We define the difference matrix $\underline{\Delta} \in \mathbb{R}^{N \times N+1}$ as
\begin{equation}
    \label{eq:delta_matrix}
    \underline{\Delta} = 
    \begin{pmatrix}
        -1 & 1 & 0 & 0 & \dots \\
        0 & -1 & 1 & 0 & \dots \\
        0 & 0 & -1 & 1 & \dots \\
        \vdots & \vdots & \vdots & \vdots & \ddots \\
    \end{pmatrix}.
\end{equation}
This difference matrix is a derivative matrix from the coefficients in the Lagrange basis to the coefficients in the histopolation basis encoding the derivative formula \eqref{eq:derivative_formula}. Indeed, for a coefficient vector $\underline{f} \vcentcolon = \left(f_0,\dots,f_N\right)^T$ of a function $f \in U_l$ in the Lagrange basis we obtain
\begin{equation}
    \underline{\Delta} \underline{f} = \left(f_1 - f_0, \dots, f_N-f_{N-1}\right)^T,
\end{equation}
where $\left(f_1 - f_0, \dots, f_N-f_{N-1}\right)^T$ is the coefficient vector of the derivative  of $f$ in the histopolation basis. 
The difference matrix together with reduction operators $R^0$ and $R^1$,
\begin{align}
    R^0: U_l \longrightarrow \mathbb{R}^{N+1}; \qquad \sum_{i = 0}^N f_i l_i(x) \longmapsto \left(f_0, \dots, f_N\right)^T, \\
    R^1: U_h \longrightarrow \mathbb{R}^{N}; \qquad \sum_{i = 1}^N \overline{g}_i h_i(x) \longmapsto \left(\overline{g}_1, \dots, \overline{g}_N\right)^T,
\end{align}
form the commuting diagram of Fig. \ref{fig:1d_discrete_de_rham}, which can be shown using Eq. \eqref{eq:derivative_formula}, relating the discrete functions with the operation on the DOFs in the discrete de Rham complex in 1D.
With Eq. \eqref{eq:derivative_formula} evaluated at a fixed grid point we also obtain
\begin{equation}
\begin{split}
    \left(\underline{D} \textbf{f}\right)_i &= \sum_{k = 0}^{N} D_{i,k} f_k  = -\sum_{j=1}^N \left(\left(f_{j}-f_{j-1}\right)\sum_{k = 0}^{j-1}D_{i,k}\right) \\
    &= \sum_{j=1}^N \left(\sum_{k = 0}^N \Delta_{j,k} f_k\right)V_{i,j} = \left(\underline{V}\, \underline{\Delta} \textbf{f}\right)_i,
\end{split}
\end{equation}
by substituting the matrix entries, where the indices of the rectangular matrix $\underline{\Delta}$ run from $1$ to $N$ for the rows and from $0$ to $N$ for the columns.
This leads to the important identity
\begin{equation}
    \label{eq:derivative_vandermonde_equality}
    \underline{D} = \underline{V}\, \underline{\Delta},
\end{equation}
which relates the nodal derivative matrix $\underline{D}$ with the histopolation Vandermonde matrix $\underline{V}$ and the difference matrix $\underline{\Delta}$.
This identity further relates the Lagrange coefficients to the histopolation coefficients of the derivative. We can now use these representations of the derivative to describe derivative operators in higher dimensions.

\subsection{Derivative Operators}

The main idea of what follows is motivated by the work of Gerritsma \cite{gerritsma2011}, and Gerritsma et al. \cite{preprint_gerritsma2011} which we aim to combine with our findings on the discrete relationships of $\underline{D}$, $\underline{V}$, and $\underline{\Delta}$ for general Lagrange and histopolation basis functions, which do not have to be polynomials.
In particular, we are interested in high-order SBP-FD operators, for which we know the basis function values and their derivatives at the grid points. In the following, vector-valued derivative operators will be denoted by bold fonts.

To represent a function $\phi$ in 3D, we use a tensor product ansatz, using the space $U_l \otimes U_l \otimes U_l$,
\begin{equation}
    \label{eq:3D_Lagrange_polynomial}
    \phi(x,y,z) = \sum_{i = 0}^{N} \sum_{j = 0}^{N} \sum_{k = 0}^{N} \phi_{i,j,k} l_i(x) l_j(y) l_k(z).
\end{equation}

Considering partial derivatives of $\phi$, we obtain, for instance, for the derivative in the $x-$direction
\begin{equation}
        \label{eq:histopolation_deriv}
        \bigpartialderiv{\phi}{x}(x,y,z) = \sum_{i = 1}^{N} \sum_{j = 0}^{N} \sum_{k = 0}^{N} (\phi_{i,j,k} - \phi_{i-1,j,k}) h_i(x) l_j(y) l_k(z).
\end{equation}
with the other derivatives evaluated analogously. 

As described in \cite{gerritsma2011}, the histopolation derivative \eqref{eq:histopolation_deriv} can be used directly to define a discrete \textbf{grad} operator from the space $\left(U_l \otimes U_l \otimes U_l\right)^T$ to the space $(U_h \otimes U_l \otimes U_l, U_l \otimes U_h \otimes U_l, U_l \otimes U_l \otimes U_h)$. 
The approximation spaces for the \textbf{curl} and div operators are adjusted accordingly with tensor product combinations of Lagrange and histopolation functions, where the general rule is that the directions in which derivatives are applied are represented by $U_l$ and other directions are represented by $U_h$.
For a more detailed explanation, see \cite{preprint_gerritsma2011}. 

As mentioned, our goal is to obtain purely discrete versions of these approximations at grid nodes for our mimetic SBP-FD approximations. It is straight forward to use the relationships in section \ref{subsection:basis_functions} to obtain analogous expressions for \textbf{grad}, \textbf{curl}, and div that are purely nodal-based. 
We do note however, that the underlying degrees of freedom are in general not only nodal, but also consist of line, area and volume integral degrees of freedom of the cell complex induced by the tensor product grid.
The derivative matrix however gives us a way of computing nodal values from the integral degrees of freedom.

For instance, the partial derivative in $x-$direction of the function $\phi$ evaluated at the grid points can be computed as
\begin{equation}
    \bigpartialderiv{}{x}\phi(x_i,y_j,z_k) = \sum_{s = 1}^{N} V_{i,s} \left( \phi_{s,j,k} - \phi_{s-1,j,k} \right) = \sum_{s = 0}^{N} D_{i,s} \phi_{s,j,k}.
\end{equation}
Note that we can compute these operators without knowing the underlying basis functions. 

In the following, we present discrete versions of all operators described by Gerritsma \cite{gerritsma2011}, which can be computed using arbitrary Lagrange and histopolation basis functions, for which we only need to know the function values and their derivatives at the grid nodes. On the other hand, as long as the values of the basis functions at the grid nodes are known, we can build the operators without having access to explicit basis functions. 
For instance, we can use high-order SBP-FD operators.

\subsubsection{Gradient Operator}

Assume the function $\phi$ from the tensor product space $U_l \otimes U_l \otimes U_l$ given through the nodal grid values $\{\phi_{i,j,k}\}_{i,j,k=0}^N$. 
Let $\textbf{grad}(\phi) = \textbf{u} = \left(u^x, u^y, u^z\right)^T$, where  $u^x \in U_h \otimes U_l \otimes U_l$, $u^y \in U_l \otimes U_h \otimes U_l$, and $u^z \in U_l \otimes U_l \otimes U_h$. Then the nodal values of the components of $\textbf{u}$ are
\begin{align}
    \label{eq:discrete_grad_x}
    u^x(x_i,y_j,z_k) = &\sum_{s = 1}^{N} V_{i,s} \left( \phi_{s,j,k} - \phi_{s-1,j,k} \right) = \sum_{s = 0}^{N} D_{i,s} \phi_{s,j,k}, \\
    \label{eq:discrete_grad_y}
    u^y(x_i,y_j,z_k) = &\sum_{s = 1}^{N} {V_{j,s}} \left( \phi_{i,s,k} - \phi_{i,s-1,k} \right) = \sum_{s = 0}^{N} D_{j,s} \phi_{i,s,k},\\
    \label{eq:discrete_grad_z}
    u^z(x_i,y_j,z_k) = &\sum_{s = 1}^{N} {V_{k,s}} \left( \phi_{i,j,s} - \phi_{i,j,s-1} \right) = \sum_{s = 0}^{N} D_{k,s} \phi_{i,j,s}.
\end{align}
\subsubsection{Curl Operator}

We aim to define the \textbf{curl} operator in a space compatible with the gradient operator described above, so that we preserve the well-known identity that the curl of the gradient of a scalar field is zero.
We choose $\tilde{\textbf{u}} = (\tilde{u}^x,\tilde{u}^y,\tilde{u}^z)$ from the same space as $\textbf{grad}(\phi)$ defined above.
That is, $\tilde{u}^x$ is from the tensor product space $U_h \otimes U_l \otimes U_l$, $\tilde{u}^y$ from $U_l \otimes U_h \otimes U_l$, and $\tilde{u}^z$ from $U_l \otimes U_l \otimes U_h$.
Then the nodal values of $\textbf{curl}(\tilde{\textbf{u}}) = \textbf{w} = (w^x,w^y,w^z)^T$ are given by
\begin{align}
    \begin{split}
    w^x(x_i,y_j,z_k) &= \sum_{s = 1}^{N} \sum_{t = 1}^{N} \left(\tilde{u}_{i,s,t}^{z} - \tilde{u}_{i,s-1,t}^{z} - \tilde{u}_{i,s,t}^{y} + \tilde{u}_{i,s,t-1}^{y}\right) V_{j,s} V_{k,t} \\
    &= \sum_{s = 0}^{N} \sum_{t = 1}^{N} \tilde{u}_{i,s,t}^{z} D_{j,s} V_{k,t} - \sum_{s = 1}^{N} \sum_{t = 0}^{N} \tilde{u}_{i,s,t}^{y} V_{j,s} D_{k,t},
    \end{split}\\
    \begin{split}
    w^y(x_i,y_j,z_k) &= \sum_{s = 1}^{N} \sum_{t = 1}^{N} \left(\tilde{u}_{s,j,t}^{x} - \tilde{u}_{s,j,t-1}^{x} - \tilde{u}_{s,j,t}^{z} + \tilde{u}_{s-1,j,t}^{z}\right) V_{i,s} V_{k,t} \\
    &= \sum_{s = 1}^{N} \sum_{t = 0}^{N} \tilde{u}_{s,j,t}^{x} V_{i,s} D_{k,t} - \sum_{s = 0}^{N} \sum_{t = 1}^{N} \tilde{u}_{s,j,t}^{z} D_{i,s} V_{k,t},
    \end{split}\\
    \begin{split}
    w^z(x_i,y_j,z_k) &= \sum_{s = 1}^{N} \sum_{t = 1}^{N} \left(\tilde{u}_{s,t,k}^{y} - \tilde{u}_{s-1,t,k}^{y} - \tilde{u}_{s,t,k}^{x} + \tilde{u}_{s,t-1,k}^{x}\right) V_{i,s}V_{j,t} \\
    &= \sum_{s = 0}^{N} \sum_{t = 1}^{N} \tilde{u}_{s,t,k}^{y} D_{i,s} V_{j,t} - \sum_{s = 1}^{N} \sum_{t = 0}^{N} \tilde{u}_{s,t,k}^{x} V_{i,s} D_{j,t}.
    \end{split}
\end{align}
If $\tilde{\textbf{u}}$ is indeed computed as the gradient of a scalar field, we can insert the coefficients from the gradient operator in Eqs. \eqref{eq:discrete_grad_x}, \eqref{eq:discrete_grad_y} and \eqref{eq:discrete_grad_z} for the coefficients of $\tilde{\textbf{u}}$.
We obtain for the coefficient of the first component of the curl
\begin{equation}
    \begin{split}
        &u_{i,j,k}^z - u_{i,j-1,k}^z - u_{i,j,k}^y + u_{i,j,k-1}^y \\
        = &\phi_{i,j,k} - \phi_{i,j,k-1} - \phi_{i,j-1,k} + \phi_{i,j-1,k-1}\\
        - &\phi_{i,j,k} + \phi_{i,j-1,k} + \phi_{i,j,k-1} - \phi_{i,j-1,k-1}\\
        = &0,
    \end{split}
\end{equation}
for all indices $i \in \{0,\dots,N\}$ and $j,k \in \{1,\dots,N\}$.
Analogously the coefficients cancel out to zero for the other two components of the curl. \\
As a consequence we have for our discrete FD operators
\begin{equation}
    \textbf{curl}(\textbf{u}) = \textbf{curl}(\textbf{grad}(\phi)) = \textbf{0},
\end{equation}
or, respectively 
\begin{equation}
    \textbf{curl} \circ \textbf{grad} = \textbf{0}.
\end{equation}

\subsubsection{Divergence Operator}

Again, we want to calculate div in a space that is compatible with the above \textbf{curl} operator, so that we preserve the well-known identity that the divergence of the curl of a vector field is zero.
We choose $\tilde{\textbf{w}} = (\tilde{w}^x,\tilde{w}^y,\tilde{w}^z)$ from the same space as $\textbf{curl}(\tilde{\textbf{u}})$ defined above.
I.e., $\tilde{w}^x$ is from the tensor product space $U_l \otimes U_h \otimes U_h$, $\tilde{w}^y$ from $U_h \otimes U_l \otimes U_h$, and $\tilde{w}^z$ from $U_h \otimes U_h \otimes U_l$. Then the nodal values of $\text{div}(\tilde{\textbf{w}}) = f$ are given by
\begin{equation}
    \label{eq:discrete_3d_div}
    \begin{split}
        f(x_i,y_j,z_k) = &\sum_{s = 1}^N \sum_{t = 1}^N \sum_{r = 1}^N (\tilde{w}_{s,t,r}^{x} - \tilde{w}_{s-1,t,r}^{x} + \tilde{w}_{s,t,r}^{y} \\
        &- \tilde{w}_{s,t-1,r}^{y} + \tilde{w}_{s,t,r}^{z} - \tilde{w}_{s,t,r-1}^{z}) V_{i,s}V_{j,t}V_{k,r} \\
        = &\sum_{s = 0}^N \sum_{t = 1}^N \sum_{r = 1}^N \tilde{w}_{s,t,r}^{x}D_{i,s}V_{j,t}V_{k,r} \\
        + &\sum_{s = 1}^N \sum_{t = 0}^N \sum_{r = 1}^N\tilde{w}_{s,t,r}^{y} V_{i,s}D_{j,t}V_{k,r}\\
        + &\sum_{s = 1}^N \sum_{t = 1}^N \sum_{r = 0}^N\tilde{w}_{s,t,r}^{z} V_{i,s}V_{j,t}D_{k,r}.
    \end{split}
\end{equation}
Analogously to the $\textbf{curl}$ operator, we can insert the coefficients of the curl into the divergence operator to obtain
\begin{equation}
   \text{div}(\textbf{w}) = \text{div}(\textbf{curl}(\tilde{\textbf{u}})) = \textbf{0},
\end{equation}
or, respectively
\begin{equation}
    \text{div} \circ \textbf{curl} = \textbf{0}.
\end{equation}

\subsubsection{2D Operators}

\begin{figure}
    \centering
    \begin{tikzpicture}
    	\node (1) at (2,0) {$H^1$};
    	\node (2) at (4,0) {$H\left(\text{div}\right)$};
    	\node (3) at (6,0) {$L^2$};
    	\draw[->] (1) -- node [midway, above] {\textbf{curl}} (2);
    	\draw[->] (2) -- node [midway, above] {div} (3);
    \end{tikzpicture}
    \caption{2D de Rham complex}
    \label{fig:2d_de_rham}
\end{figure}
\begin{figure}
    \centering
    \begin{tikzpicture}
    	\node (1) at (2,0) {$L^2$};
    	\node (2) at (4,0) {$H\left(\text{curl}\right)$};
    	\node (3) at (6,0) {$H^1$};
    	\draw[<-] (1) -- node [midway, above] {curl} (2);
    	\draw[<-] (2) -- node [midway, above] {-\textbf{grad}} (3);
    \end{tikzpicture}
    \caption{2D dual de Rham complex}
    \label{fig:2d_dual_de_rham}
\end{figure}
In 2D, we have four relevant derivative operators on our discrete spaces that we can use to find a discrete counterpart to the 2D primal and dual de Rham complexes of Fig.~\ref{fig:2d_de_rham} and Fig.~\ref{fig:2d_dual_de_rham}. The 2D $\textbf{grad}$ and $\text{div}$ operators, one scalar-to-vector $\textbf{curl} = \left(\bigpartialderiv{}{y}, -\bigpartialderiv{}{x}\right)^T$ and a vector-to-scalar $\text{curl}\left(\left(v^x, v^y\right)^T\right) = \bigpartialderiv{v^y}{x} - \bigpartialderiv{v^x}{y}$. 
Both $\textbf{grad}$ and $\textbf{curl}$ are defined for functions $u \in U_l \otimes U_l$ while div is defined for functions $\textbf{v} \in \left(U_l \otimes U_h, U_h \otimes U_l \right)^T$.
Finally, the curl operator is defined for functions $\textbf{w} \in \left(U_h \otimes U_l, U_l \otimes U_h \right)^T$.

In all cases, we have Lagrange basis functions in the directions in which we take the derivatives and histopolation functions otherwise.
In the following, we will write the degrees of freedom for all functions in vector form, where we save them in the following order:
\begin{figure}
    \centering
    \resizebox{\textwidth}{!}{
    \begin{tikzpicture}
        \def\N{4}
        \pgfmathsetmacro{\NmOne}{\N-1}
        \def\spacing{6}  

        \begin{scope}[xshift=0cm]  
            \draw[thick] (0,0) rectangle (\N, \N);

            \foreach \x in {0,...,\N} {
                \foreach \y in {0,...,\N} {
                    \draw[black] (\x, \y) +(-0.1,-0.1) rectangle +(0.1,0.1);
                }
            }

            \foreach \x in {0,...,\N} {
                \foreach \y in {0,...,\NmOne} {
                    \fill[black] (\x-0.05, \y+0.15) rectangle (\x+0.05, \y+0.85);
                }
            }

            \node[anchor=east] at (0,0.5) {$(0,1)$};
            \node[anchor=east] at (1,0.5) {$(1,1)$};
            \node[anchor=east] at (0,1.5) {$(0,2)$};
            \node[anchor=west] at (\N,0.5) {$(N,1)$};
            \node[anchor=west] at (\N,\N-0.5) {$(N,N)$};
            \node[anchor=east] at (0,\N-0.5) {$(0,N)$};
            \node[anchor=east] at (2,0.5) {$\cdots$};
            \node[anchor=east] at (-0.4,2.5) {$\vdots$};
            \node[anchor=west] at (\N+0.4,1.5) {$\vdots$};

            \node[anchor=north] at (\N/2, -0.5) {$v^x$};
        \end{scope}

        \begin{scope}[xshift=\spacing cm] 
            \draw[thick] (0,0) rectangle (\N, \N);

            \foreach \x in {0,...,\N} {
                \foreach \y in {0,...,\N} {
                    \draw[black] (\x, \y) +(-0.1,-0.1) rectangle +(0.1,0.1);
                }
            }

            \foreach \x in {0,...,\NmOne} {
                \foreach \y in {0,...,\N} {
                    \fill[black] (\x+0.15, \y-0.05) rectangle (\x+0.85, \y+0.05);
                }
            }

            \node[anchor=north] at (0.5,0) {$(1,0)$};
            \node[anchor=north] at (1.5,0) {$(2,0)$};
            \node[anchor=north] at (0.5,1) {$(1,1)$};
            \node[anchor=south] at (0.5,\N) {$(1,N)$};
            \node[anchor=south] at (\N-0.5,\N) {$(N,N)$};
            \node[anchor=north] at (\N-0.5,0) {$(N,0)$};
            \node[anchor=north] at (2.5,0) {$\cdots$};
            \node[anchor=north] at (0.5,2.0) {$\vdots$};
            \node[anchor=south] at (1.5,\N) {$\cdots$};

            \node[anchor=north] at (\N/2, -0.5) {$v^y$};
        \end{scope}

        \begin{scope}[xshift=2*\spacing cm]
            \draw[thick] (0,0) rectangle (\N, \N);

            \foreach \x in {0,...,\N} {
                \foreach \y in {0,...,\N} {
                    \draw[black] (\x, \y) +(-0.1,-0.1) rectangle +(0.1,0.1);
                    \fill[black] (\x, \y) circle (0.05);
                }
            }

            \node[anchor=north east] at (0,0) {$(0,0)$};
            \node[anchor=east] at (0,1) {$(0,1)$};
            \node[anchor=north] at (1,0) {$(1,0)$};

            \node[anchor=east] at (-0.3,2) {$\vdots$};
            \node[anchor=north] at (2,-0.1) {$\cdots$};
            \node[anchor=south west] at (2,2) {\text{\rotatebox[origin=c]{45}{$\cdots$}}}; 

            \node[anchor=south] at (1,\N) {$\cdots$};
            \node[anchor=west] at (\N+0.2,1) {$\vdots$};

            \node[anchor=north west] at (\N,0) {$(N,0)$};
            \node[anchor=south west] at (\N,\N) {$(N,N)$};
            \node[anchor=south east] at (0,\N) {$(0,N)$};

            \node[anchor=north] at (\N/2, -0.5) {$u$};
        \end{scope}
    \end{tikzpicture}
    }
    \caption{Storage of the degrees of freedom for 2D functions}
    \label{fig:varstorage}
\end{figure}
{
\begin{align}
    \textbf{\underline{v}}^x& = 
    \begin{bmatrix}
        v^x_{01} \\
        \vdots \\
        v^x_{N1} \\
        v^x_{02} \\
        \vdots \\
        v^x_{N2} \\
        \vdots \\
        v^x_{NN}
    \end{bmatrix}
    \in \mathbb{R}^{N(N+1)},\qquad
    \textbf{\underline{v}}^y = 
    \begin{bmatrix}
        v^y_{10} \\
        \vdots \\
        v^y_{N0} \\
        v^y_{11} \\
        \vdots \\
        v^y_{N1} \\
        \vdots \\
        v^y_{NN}
    \end{bmatrix}
    \in \mathbb{R}^{N(N+1)}, \label{eq:disc_storage_mixed}\\
    \underline{u} &= 
    \begin{bmatrix}
        u_{00} \\
        \vdots \\
        u_{N0} \\
        u_{01} \\
        \vdots \\
        u_{N1} \\
        \vdots \\
        u_{NN}
    \end{bmatrix}
    \in \mathbb{R}^{(N+1)^2},
    \label{eq:disc_storage_Lagrange}
\end{align}
}
where $\textbf{v}$ is a function for which we want to be able to compute the divergence. 
Analogously to $\mathbf{v}$, we can define the discrete storage of a function $\textbf{w}$, for which we need to compute its vector-to-scalar curl.
Figure \ref{fig:varstorage} shows in which order the DOFs are indexed.
Using the Kronecker product, we can define the 2D partial derivative operators in matrix form:
\begin{align}
   \underline{\Delta}^{2D}_x &\coloneqq \underline{I}_{M} \otimes \underline{\Delta} \in {\mathbb{R}^{M N \times M (N+1)}}, \\
   \underline{\Delta}^{2D}_y &\coloneqq  \underline{\Delta} \otimes \underline{I}_{M} \in {\mathbb{R}^{M N \times M (N+1)}},
\end{align}
where $\underline{I}_{M}$ is the $M \times M$ identity matrix.
We have $M = N + 1$ if the derivative acts on $U_l \otimes U_l$ or $M = N$ if it acts on $U_l \otimes U_h$ or $U_h \otimes U_l$.
In an abuse of notation we do not add the value of $M$ as an index of the $\underline{\Delta}^{2D}$ matrices but leave it implied by context.

Now we can write the derivative operators in matrix form:
\begin{align}
    \underline{\textbf{grad}}^x\left(\underline{u}\right) &= \underline{\Delta}^{2D}_x \underline{u}, \\
    \underline{\textbf{grad}}^y\left(\underline{u}\right) &= \underline{\Delta}^{2D}_y \underline{u}, \\
    \label{eq:2d_discrete_curl_x}
    \underline{\textbf{curl}}^x\left(\underline{u}\right) &= \underline{\Delta}^{2D}_y \underline{u}, \\
    \label{eq:2d_discrete_curl_y}
    \underline{\textbf{curl}}^y\left(\underline{u}\right) &= -\underline{\Delta}^{2D}_x \underline{u}, \\
    \underline{\text{curl}}\left(\textbf{\underline{w}}^x,\textbf{\underline{w}}^y\right) &= \underline{\Delta}^{2D}_x \textbf{\underline{w}}^y - \underline{\Delta}^{2D}_y \textbf{\underline{w}}^x, \\
    \label{eq:2d_discrete_div}
    \underline{\text{div}}\left(\textbf{\underline{v}}^x,\textbf{\underline{v}}^y\right) &= \underline{\Delta}^{2D}_x \textbf{\underline{v}}^x + \underline{\Delta}^{2D}_y \textbf{\underline{v}}^y,
\end{align}
and we obtain
\begin{align}
    \label{eq:2D_curl_div_identity}
    \underline{\text{curl}} \circ \underline{\textbf{grad}}\left(\underline{u}\right) &= { \left(\underline{\Delta}^{2D}_x \underline{\Delta}^{2D}_y - \underline{\Delta}^{2D}_y \underline{\Delta}^{2D}_x\right)\underline{u}} = \underline{0},\\
    \label{eq:2D_grad_curl_identity}
    \underline{\text{div}} \circ \underline{\textbf{curl}}\left(\underline{u}\right) &= \left(\underline{\Delta}^{2D}_x \underline{\Delta}^{2D}_y - \underline{\Delta}^{2D}_y \underline{\Delta}^{2D}_x\right)\underline{u} = \underline{0},   
\end{align}
where we have used the following property of Kronecker products:
\begin{equation}
    \left(\underline{A} \otimes \underline{B}\right)\left(\underline{C} \otimes \underline{D}\right) = \underline{A}\underline{C} \otimes \underline{B}\underline{D}.
\end{equation}

With the matrix representation of our derivative operators, we have the commuting 2D de Rham complex of Fig.~\ref{fig:2d_discrete_de_rham}, using reduction operators $R^0, \textbf{R}^1$ and $R^2$.
These operators are defined as the mappings sending functions in the discrete spaces to the vector of their degrees of freedom.
We can again evaluate the operators nodally using the histopolation Vandermonde matrix $\underline{V}$ analogously to the 3D case in equations \eqref{eq:discrete_grad_x} to \eqref{eq:discrete_3d_div}.

\begin{figure}
    \centering
    \begin{tikzpicture}
        \label{derham2d_feec}
    	\node (1) at (3,0) {$H^1$};
     	\node (2) at (6.5,0) {$H(\text{div})$};
       	\node (3) at (10,0) {$L^2$};
       	\node (6) at (3,-1.5) {$W^0$};
       	\node (7) at (6.5,-1.5) {$W^1$};
       	\node (8) at (10,-1.5) {$W^2$};
       	\node (9) at (3,-3) {$\mathbb{R}^{(N+1)^2}$};
       	\node (10) at (6.5,-3) {$\mathbb{R}^{(N+1)N} \times \mathbb{R}^{N(N+1)}$};
       	\node (11) at (10,-3) {$\mathbb{R}^{N^2}$};
       	\draw[->] (1) -- node [midway, above] {\textbf{curl}} (2);
       	\draw[->] (2) -- node [midway, above] {div} (3);
       	\draw[->] (7) -- node [midway, above] {div} (8);
       	\draw[->] (6) -- node [midway, above] {\textbf{curl}} (7);
       	\draw[->] (10) -- node [midway, above] {\underline{div}} (11);
       	\draw[->] (9) -- node [midway, above] {\underline{\textbf{curl}}} (10);
       	\draw[->] (1) -- node [midway, right] {$p^0$} (6);
       	\draw[->] (3) -- node [midway, right] {$p^2$} (8);
       	\draw[->] (2) -- node [midway, right] {$\textbf{p}^1$} (7);
       	\draw[->] (6) -- node [midway, right] {$R^0$} (9);
       	\draw[->] (7) -- node [midway, right] {$\textbf{R}^1$} (10);
       	\draw[->] (8) -- node [midway, right] {$R^2$} (11);
    \end{tikzpicture}
    \caption{Discrete de Rham complex in 2D}
    \label{fig:2d_discrete_de_rham}
\end{figure}
\newpage
\begin{figure}
    \centering
    \begin{tikzpicture}
    	\node (1) at (2,0) {$H^1$};
    	\node (2) at (4,0) {$H\left(\text{curl}\right)$};
    	\node (3) at (6.5,0) {$H\left(\text{div}\right)$};
        \node (4) at (8.5,0) {$L^2$};
    	\draw[->] (1) -- node [midway, above] {\textbf{grad}} (2);
    	\draw[->] (2) -- node [midway, above] {\textbf{curl}} (3);
    	\draw[->] (3) -- node [midway, above] {div} (4);
    \end{tikzpicture}
    \caption{{3D de Rham complex}}
    \label{fig:3d_de_rham}
\end{figure}
One further object of interest are the cohomology spaces of the de Rham complexes in Fig. \ref{fig:2d_de_rham}, Fig \ref{fig:2d_dual_de_rham} and Fig. \ref{fig:3d_de_rham}.
For vector spaces $V \subset U$ we denote by $\faktor{U}{V}$ the quotient space of $U$ to its subspace $V$.
Then the cohomology spaces for the 2D complex in Fig. \ref{fig:2d_de_rham} are defined as 
\begin{align}
    \label{eq:cohomology_primal_0}
    H^0 \vcentcolon &= \faktor{\text{Ker}\left(\textbf{curl}\right)}{I\left(\mathbb{R}\right)}, \\
    \label{eq:cohomology_primal_1}
    H^1 \vcentcolon &= \faktor{\text{Ker}\left(\text{div}\right)}{\text{Im}\left(\textbf{curl}\right)}, \\
    \label{eq:cohomology_primal_2}
    H^2 \vcentcolon &= \faktor{L^2}{\text{Im}\left(\text{div}\right)},
\end{align}
where $\text{Ker}$ denotes the kernel, $\text{Im}$ denotes the image, and $I\left(\mathbb{R}\right)$ is the inclusion of $\mathbb{R}$ into $H^1$ by mapping $c \in \mathbb{R}$ to the constant function with value $c$.
Analogously we have the cohomology spaces for the dual complex in Fig. \ref{fig:2d_dual_de_rham}
\begin{align}
    \label{eq:cohomology_dual_0}
    \tilde{H}^0 \vcentcolon &= \faktor{\text{Ker}\left(\textbf{grad}\right)}{I\left(\mathbb{R}\right)}, \\
    \label{eq:cohomology_dual_1}
    \tilde{H}^1 \vcentcolon &= \faktor{\text{Ker}\left(\text{curl}\right)}{\text{Im}\left(\textbf{grad}\right)}, \\
    \label{eq:cohomology_dual_2}
    \tilde{H}^2 \vcentcolon &= \faktor{L^2}{\text{Im}\left(\text{curl}\right)},
\end{align}
and in 3D for the complex in Fig. \ref{fig:3d_de_rham}
\begin{align}
    \label{eq:cohomology_3D_0} 
    H^0_{3D} &\vcentcolon = \faktor{\text{Ker}\left(\textbf{grad}\right)}{I\left(\mathbb{R}\right)}, \\ 
    \label{eq:cohomology_3D_1}
    H^1_{3D} &\vcentcolon = \faktor{\text{Ker}\left(\textbf{curl}\right)}{\text{Im}\left(\textbf{grad}\right)}, \\
    \label{eq:cohomology_3D_2}
     H^2_{3D} &\vcentcolon = \faktor{\text{Ker}\left(\text{div}\right)}{\text{Im}\left(\textbf{curl}\right)}, \\
    \label{eq:cohomology_3D_3}
    H^3_{3D} &\vcentcolon = \faktor{L^2}{\text{Im}\left(\text{div}\right)}.
\end{align}
The dimensions of these spaces are topological invariants only depending on the domain $\Omega$.
The discrete cohomology spaces are defined analogously replacing the analytic derivative operators with their discrete analogues (e.g. $\underline{\textbf{grad}}$ instead of  $\textbf{grad}$) and replacing $L^2$ with the discrete space $W^2$ (or $W^3$ in 3D). 
The definition of $I(\mathbb{R})$ remains the same, since by construction the space $W^0$ contains the constant functions due to Eq. \eqref{eq:Lagrange_constant_functions}.
We denote these discrete cohomology spaces with the subscript $h$, e.g. $H^0_h$ is the discrete analogue to $H^0$.
Finite Element Exterior Calculus preserves the dimension of these cohomology spaces in the discrete complex when the projections fulfill fairly lenient requirements (see \cite{arnold2010}, theorem 3.4).
Since our finite difference operators are defined using a similar approach to FEEC, we are therefore interested if the discrete complex defined via our finite difference operators also preserve the dimension of the cohomology spaces.
For the case where $\Omega$ is a non-periodic rectangular or hexahedral domain, we can explicitly prove that the discrete cohomology spaces do indeed vanish like the analytic ones:
\begin{proposition}
    For a non-periodic rectangular domain all cohomology spaces of the cohomology spaces $H_h^0$, $H_h^1$, $H_h^2$, $\tilde{H}_h^0$, $\tilde{H}_h^1$, and $\tilde{H}_h^2$ vanish, and are therefore the same as for the analytic de Rham complexes.
    For a non-periodic hexahedral domain the cohomology spaces $H_{3D,h}^0$, $H_{3D,h}^1$, $H_{3D,h}^2$, and $H_{3D,h}^3$ also vanish in agreement with the analytic complex.\\
    \\
    \textbf{Proof:}\\
    The proof works by counting degrees of freedom. For details see appendix A.
\end{proposition}

\subsection{Summation-By-Parts Finite Difference Operators} \label{sec:SBP-FD}
In the previous section, we assumed the existence of Lagrange functions that fulfill {Eq. \eqref{eq:Lagrange_constant_functions} and Eq. \eqref{eq:Lagrange_linear_functions}} at grid points $\{x_i\}_{i = 0}^N$ (including the boundary nodes), but we did not discuss how to construct them.
In \cite{preprint_gerritsma2011}, Gerritsma uses the Legendre--Gauss--Lobatto points and the Lagrange {interpolation} polynomials for those points.
The corresponding derivative matrix $\underline{D}$ together with the mass matrix $\underline{M}$, which is a diagonal matrix with the Legendre--Gauss--Lobatto quadrature weights on the diagonal, has the summation-by-parts (SBP) property,
\begin{equation}
    \label{eq:SBP}
    \underline{M}\underline{D} + \underline{D}^T\underline{M} = \underline{B},
\end{equation}
where
\begin{equation}
    B_{i,j} = -\delta_{0,i}\delta_{0,j} + \delta_{N,i}\delta_{N,j}
\end{equation}
is a boundary matrix.

As we have pointed out above, we only need to evaluate our basis functions and their derivatives at the grid nodes to construct the mimetic operators defined in the previous section. 
If we use a collocation quadrature rule when applying Finite Element Exterior Calculus \cite{arnold2006,arnold2010} to derive our scheme, the global shape of the basis functions $l_i$ does not matter at all, only their nodal function and derivative values are relevant.
The central idea of this paper is therefore to not construct basis functions $l_i$ at all, but use existing diagonal norm finite difference SBP operators, for example from Strand \cite{strand1994}, to supply a quadrature rule and nodal derivative data needed for the numerical scheme.
We do this explicitly by prescribing the nodal derivatives of the Lagrange functions, that are not known globally, through the derivative matrix according to {Eq. \eqref{eq:derivative_matrix}} and using diagonal entries of the norm matrix $\underline{M}$ as weights for a collocation quadrature rule on the grid points.

One of the standard SBP finite difference operators from \cite{strand1994}, which we will use in this study, is the operator with interior order $4$ and boundary order $2$ with $N+1$ points:
\begin{align}
    \label{eq:SBP_derivative_matrix}
    \underline{D} &= \frac{1}{\Delta x}\begin{pmatrix}-\frac{24}{17} & \frac{59}{34} & \frac{4}{17} & -\frac{3}{34} & 0 & 0 & 0 & \dots \\
    -\frac{1}{2} & 0 & \frac{1}{2} & 0 & 0 & 0 & 0 & \dots \\
    \frac{4}{43} & -\frac{59}{86} & 0 & \frac{59}{86} & -\frac{4}{43} & 0 & 0 & \dots\\
    \frac{3}{98} & 0 & -\frac{59}{98} & 0 & \frac{32}{49} & -\frac{4}{49} & 0 & \dots \\
    0 & 0 & \frac{1}{12} & -\frac{2}{3} & 0 & \frac{2}{3} & -\frac{1}{12} & \dots \\
     & & & \ddots & \ddots & \ddots & \ddots & \ddots
    \end{pmatrix},\\
    \label{eq:SBP_mass_matrix}
    \underline{M} &= \Delta x \hspace{0.05cm} \text{diag}\left(\frac{17}{48}, \frac{59}{48}, \frac{43}{48}, \frac{49}{48}, 1, \dots, 1, \frac{49}{48}, \frac{43}{48}, \frac{59}{48}, \frac{17}{48}\right),
\end{align}
which is defined on a regular grid, where $\Delta x$ is the distance between two neighboring grid nodes.
In the derivative matrix, the upper left boundary block is repeated in the bottom right with reversed signs and inverted order.

In the following, we will derive a mimetic scheme for the 2D transverse electric Maxwell's equations to validate our approach.

\section{Mimetic Scheme for the 2D Maxwell's Equations}

\subsection{Maxwell's Equations in 2D}

The homogeneous transverse electric Maxwell's equations with unit light speed read
\begin{align}
    \label{eq:faraday}
    \bigpartialderiv{B^z}{t} + \text{curl}\left(\textbf{E}\right) &= 0,\\
    \label{eq:ampere}
    \bigpartialderiv{\textbf{E}}{t} - \textbf{curl}\left(B^z\right) &= \textbf{0},\\
    \label{eq:gauss_law}
    \text{div}\left(\textbf{E}\right) &= 0.
\end{align}
Our goal is to derive a mimetic scheme using formulations analogous to finite element exterior calculus to preserve the divergence involution \eqref{eq:gauss_law}.
For simplicity we will consider a rectangular domain, using a Cartesian mesh for the discretization.
To apply the framework of Finite Element Exterior Calculus, we need the de Rham complexes containing the derivative operators of the PDEs.

\subsection{De Rham Complexes}

\begin{figure}
    \centering
    \begin{tikzpicture}
    	\node (1) at (2,0) {$H_0^1$};
    	\node (2) at (4,0) {$H_0\left(\text{div}\right)$};
    	\node (3) at (6,0) {$L^2$};
    	\draw[->] (1) -- node [midway, above] {\textbf{curl}} (2);
    	\draw[->] (2) -- node [midway, above] {div} (3);
    \end{tikzpicture}
    \caption{2D de Rham complex with homogeneous boundary conditions}
    \label{fig:2d_de_rham_bc}
\end{figure}
\begin{figure}
    \centering
    \begin{tikzpicture}
    	\node (1) at (2,0) {$L^2$};
    	\node (2) at (4,0) {$H_0\left(\text{curl}\right)$};
    	\node (3) at (6,0) {$H_0^1$};
    	\draw[<-] (1) -- node [midway, above] {curl} (2);
    	\draw[<-] (2) -- node [midway, above] {-\textbf{grad}} (3);
    \end{tikzpicture}
    \caption{2D dual de Rham complex with homogeneous boundary conditions}
    \label{fig:2d_dual_de_rham_bc}
\end{figure}

Maxwell's equations have an inherent geometric structure described by the de Rham cochain complex of Fig. \ref{fig:2d_de_rham} and its dual in the periodic case in Fig. \ref{fig:2d_dual_de_rham} in 2D. We define the spaces needed for the complexes by
\begin{align}
    H\left(\text{curl}\right) &= \{\textbf{w} \in L^2 \hspace{0.05cm}|  \hspace{0.05cm} \text{curl}\left(\textbf{w}\right) \in L^2\},\\
    H\left(\text{div}\right) &= \{\textbf{w} \in L^2 \hspace{0.05cm}|  \hspace{0.05cm} \text{div}\left(\textbf{w}\right) \in L^2\}.
\end{align}
The complexes are dual for periodic domains in the sense that
\begin{align}
    \label{eq:duality_curl}
    \left< \textbf{curl}\left(\psi\right), \boldsymbol{\Phi} \right>
     &= \left< \psi, \text{curl}\left(\boldsymbol{\Phi}\right) \right> ,&\qquad \forall \psi \in H^1, \boldsymbol{\Phi} \in H\left(\text{curl}\right),\\
    \label{eq:duality_grad}
    \left< \text{div}\left(\boldsymbol{\Psi}\right), \phi \right> &= \left< \boldsymbol{\Psi}, -\textbf{grad}\left(\phi\right) \right>, &\qquad \forall \boldsymbol{\Psi} \in  H\left(\text{div}\right), \phi \in  H^1,
\end{align}
where the brackets denote the $L^2$ inner product. 
In the non-periodic case, the duality equations \eqref{eq:duality_curl} and \eqref{eq:duality_grad} do not hold any longer due to appearing boundary integrals when integrating by parts.
However, we can restore duality when considering homogeneous boundary conditions for the primal or the dual complex.
We consider the spaces
\begin{align}
    H^1_0 &= \{w \in H^1 \hspace{0.05cm}|  \hspace{0.05cm} u_{|\partial \Omega}\equiv 0\},\\
    H_0 \left(\text{curl}\right) &= \{\textbf{w} \in  H\left(\text{curl}\right) \hspace{0.05cm}|  \hspace{0.05cm} \textbf{w}_{|\partial \Omega} \cdot \tau \equiv 0 \text{ for all tangent vector fields } \tau \},\\
    H_0 \left(\text{div}\right) &= \{\textbf{w} \in  H\left(\text{div}\right) \hspace{0.05cm}|  \hspace{0.05cm} \textbf{w}_{|\partial \Omega} \cdot \textbf{n} \equiv 0 \text{ for all normal vector fields } \textbf{n}\},
\end{align}
where we again omit the region $\Omega$ in the space names for brevity.
Note that there is no adjustment for the $L^2$ space, since the trace operator is not well-defined on this space.
Imposing homogeneous boundary conditions on either the primal or the dual complex restores duality. 
In other words, the primal complex without imposed boundary conditions in Fig. \ref{fig:2d_de_rham} and the dual complex with homogeneous boundary conditions in Fig. \ref{fig:2d_dual_de_rham_bc} are dual to each other.
So are complexes in Fig. \ref{fig:2d_de_rham_bc} and in Fig. \ref{fig:2d_dual_de_rham}.
If we now demand in equations \eqref{eq:duality_curl} and \eqref{eq:duality_grad} that $\psi \in H^1_0$ and $\boldsymbol{\Psi} \in H_0 \left(\text{div}\right)$ or $\boldsymbol{\Phi} \in H_0\left(\text{curl}\right)$ and $\phi \in H^1_0$, the equations hold for non-periodic domains.
The homogeneous boundary conditions for the primal complex are:
\begin{align}
    B^z_{|\partial \Omega} &\equiv 0, \\
    \textbf{E}_{|\partial \Omega} \cdot \textbf{n} &\equiv 0,
\end{align}
for all normal vector fields $\textbf{n}$, while for the dual complex they are
\begin{equation}
    \textbf{E}_{|\partial \Omega} \cdot {\tau} \equiv 0,
\end{equation}
for all tangent vector fields {$\tau$}.
We can also impose homogeneous boundary conditions on the primal complex for one part of the boundary and on the dual complex on the remaining part.
We have to choose in which spaces to look for $\textbf{E}$ and $B^z$, since we have to choose which complex to discretize.
If we choose $\textbf{E} \in H(\text{div})$ or $\textbf{E} \in H_0(\text{div})$ and $B^z \in H^1$ or $B^z \in H_0^1$, we can solve Ampère's law and Gauss's law strongly, while we have to use Stokes' theorem to solve Faraday's law weakly.
Analogously, if we choose $\textbf{E} \in H(\text{curl})$ or $\textbf{E} \in H_0(\text{curl})$ and $B^z \in L^2$, we can solve Faraday's law strongly, but need to solve Ampère's law and Gauss's law weakly.\\
Since the divergence involution for the magnetic field falls away in the 2D transverse electric case and we want to preserve Gauss' law for the electric field \eqref{eq:gauss_law}, we will choose $\textbf{E} \in H(\text{div})$ and $B^z \in H^1$, leading to the formulation for the periodic case:
\begin{align}
    \label{eq:faraday_variational}
    \left<\bigpartialderiv{B^z}{t},\psi\right> + \left< \textbf{E}, \textbf{curl}\left( \psi \right)\right> &= 0, \qquad & \forall \psi \in H^1, \\
    \label{eq:ampere_variational}
    \left<\bigpartialderiv{\textbf{E}}{t}, \boldsymbol{\Phi}\right> - \left< \textbf{curl}\left(B^z\right), \boldsymbol{\Phi} \right> &= 0, \qquad & \forall \boldsymbol{\Phi} \in H(\text{div}),\\
    \label{eq:gauss_law_variational}
    \left<\text{div}\left(\textbf{E}\right), \mu \right> &= 0, \qquad & \forall \mu \in L^2.
\end{align}
For non-periodic problems, the function spaces must be modified according to the imposed boundary conditions. 
Boundary conditions imposed directly through the definition of the trial and test spaces correspond to \textit{essential} boundary conditions. 
In contrast, boundary conditions arising through the boundary terms generated by integration by parts correspond to \textit{natural} boundary conditions.

The finite element exterior calculus ansatz is to choose a discrete subcomplex of the de Rham complex, meaning we choose $W^0 \subset H^1$, $W^1 \subset H\left(\text{curl}\right)$, $W^2 \subset L^2$ such that $\textbf{curl}\left(W^0\right) \subset W^1$ and $\text{div}\left(W^1\right) \subset W^2$.
If we use essential boundary conditions, we instead choose subspaces of the spaces of the complex with homogeneous boundary conditions.
In addition we need bounded projections $p_0, \textbf{p}_1$ and $p_2$, so that the diagram of Fig.~\ref{fig:2d_discrete_de_rham} commutes.
We then restrict both our ansatz and our test spaces to the corresponding finite dimensional subspaces.
This means we again have two ansatz spaces for the electric and magnetic field: We look for a numerical solution $\textbf{E}_h \in [0,T] \times W^1$ and $B^z_h \in [0,T] \times W^0$ so that
\begin{align}
    \label{eq:faraday_discrete}
    \left<\bigpartialderiv{B^z_h}{t},\psi_h\right> + \left< \textbf{E}_h, \textbf{curl}\left( \psi_h \right)\right> &= 0, \qquad & \forall \psi_h \in W^0, \\
    \label{eq:ampere_discrete}
    \left<\bigpartialderiv{\textbf{E}_h}{t}, \boldsymbol{\Phi}_h\right> - \left< \textbf{curl}\left(B^z_h\right), \boldsymbol{\Phi}_h \right> &= 0, \qquad & \forall \boldsymbol{\Phi}_h \in W^1,\\
    \label{eq:gauss_discrete}
    \left<\text{div}\left(\textbf{E}_h\right), \mu_h \right> &= 0, \qquad & \forall \mu_h \in W^2.
\end{align}
Since $\text{div}(\textbf{E}_h) \in W^2$ we do obtain $\text{div}(\textbf{E}_h) \equiv 0$ from \eqref{eq:gauss_discrete}.
We will later show that by solving \eqref{eq:faraday_discrete} and \eqref{eq:ampere_discrete} with this method, we automatically preserve the divergence involution both in the semi-discrete and fully discrete setting.
Therefore we will only consider the evolution equations \eqref{eq:faraday_discrete} and \eqref{eq:ampere_discrete} when deriving our method.

\subsection{Spaces and Projections}
To construct our spaces and projections, we first look at the one-dimensional case. We have an interval $[a,b]$ and a partition of $\{K^k\}_{k = 1}^m$ of ${a,b}$ into elements.
The global $U_l$ and $U_h$ spaces are defined via the local $U_l$ and $U_h$ spaces on each element.
In the case of $U_l$ we impose global continuity to guarantee $W^0 \subset H^1$ and $W^1 \subset H(\text{div})$ (see for example Monk \cite{monk_2003}), while there are no such requirements for $U_h$.
Further details can be found in appendix C.
To extend this approach to 2D, we use a tensor product ansatz, where for a rectangular periodic domain $\Omega$ with a Cartesian grid, we define the discrete spaces via a tensor product ansatz as follows:
\begin{align}
    W^0 &:= U_l \otimes U_l, \\
    W^1 &:= \left(U_l \otimes U_h, U_h \otimes U_l\right)^T, \\
    W^2 &:= U_h \otimes U_h.
\end{align}
Here we have a slight abuse of notation because we can have different number of elements and scaling factors in both directions. 
The projections $p^0,\mathbf{p}^1,$ and $p^2$  are now defined by taking the interpolation in the directions where we use the space $U_l$ and histopolation in the directions where we have $U_h$. For the periodic case we have:
\begin{align}
    p^0(f) &:= \sum_{k = 1}^m \sum_{i,j = 0}^{N-1} f(x_i^k,y_j^k) l_i(x) l_j(y), \\
    p_1^1(\textbf{f}) &:= \sum_{k = 1}^m \sum_{i = 0}^{N-1} \sum_{j = 1}^{N} \left(\int_{y_{j-1}^k}^{y_j^k}
    f(x_i^k,t)\text{ d}t\right) l_i^k(x) h_j^k(y), \\
    p_2^1(\textbf{f}) &:= \sum_{k = 1}^m \sum_{i = 1}^{N} \sum_{j = 0}^{N-1} \left(\int_{x_{i-1}^k}^{x_i^k}
    f(s,y_j^k)\text{ d}s\right) h_i^k(x) l_j^k(y), \\
    p^2(f) &:= \sum_{k = 1}^m \sum_{i,j = 1}^{N} \left(\int_{x_{i-1}^k}^{x_i^k} \int_{y_{j-1}^k}^{y_j^k} f(s,t) \text{ d}s \text{ d}t\right) h_i^k(x) h_j^k(y),
\end{align}
where $m$ is the number of elements.
The sum index for the Lagrange functions only goes up to $N-1$ to prevent adding a basis function twice.
In the non-periodic case we may need to add or subtract Lagrange basis functions at the boundary depending on imposed boundary conditions.
Since we use the same DOFs as Gerritsma et al. \cite{preprint_gerritsma2011}, the boundedness of the projections carries over, since we can factorize over Gerritsma's projections switching out the basis functions.
As a linear mapping between finite dimensional vector spaces, it is bounded, and therefore our projections are bounded as the composition of two bounded mappings.
It is straightforward to see that we also get commuting property for the discrete 2D de Rham complex in Fig. \ref{fig:2d_discrete_de_rham}.

\section{Numerical Method}

\subsection{Galerkin Discretization}
In this and the next section we will follow standard continuous Galerkin finite element procedures to derive our scheme.
We will then use collocation quadrature to reduce the formulation to nodal and integral degrees of freedom.
Before we derive the method, we have to account for the overlap in degrees of freedom associated with the basis functions associated with the boundary nodes in continuous directions. 
To avoid duplication that we would have to account {for} in our formulation, we define the vectors of degrees of freedom as:
\begin{equation*}
    \underline{B}^z = \begin{pmatrix}
        B^z_{0,0;1} \\
        \vdots \\
        B^z_{N-1,0;1} \\
        B^z_{0,1;1} \\
        \vdots \\
        B^z_{N-1,N-1;1} \\
        B^z_{0,0;2} \\
        \vdots \\
        B^z_{N-1,N-1;m}
    \end{pmatrix},
    \qquad
    \overline{\underline{E}}^x = \begin{pmatrix}
        \overline{E}^x_{0,1;1} \\
        \vdots \\
        \overline{E}^x_{N-1,1;1} \\
        \overline{E}^x_{0,2;1} \\
        \vdots \\
        \overline{E}^x_{N-1,N;1} \\
        \overline{E}^x_{0,1;2} \\
        \vdots \\
        \overline{E}^x_{N-1,N;m}
    \end{pmatrix},
    \qquad
    \overline{\underline{E}}^x = \begin{pmatrix}
        \overline{E}^y_{1,0;1} \\
        \vdots \\
        \overline{E}^y_{N,0;1} \\
        \overline{E}^y_{1,1;1} \\
        \vdots \\
        \overline{E}^y_{N,N-1;1} \\
        \overline{E}^y_{1,0;2} \\
        \vdots \\
        \overline{E}^y_{N,N-1;m}
    \end{pmatrix}.
\end{equation*}
For ease of notation, we define $B_{N,j,k} := B_{0,j,\sigma_r(k)}$, $B_{i,N,k} := B_{i,0,\sigma_r(k)}$, and $B_{N,N,k} := B_{0,0,\sigma_r\left(\sigma_u(k)\right)}$. 
In the non-periodic case, the vector segments for the right and upper boundary elements contain the boundary value as well.
Here, $\sigma_r(k)$ is defined as the index of the right neighbor element of the element $K_k$ and $\sigma_u(k)$ is the index of the upper neighbor of $K_k$.
We also apply this notation for the electric field coefficients in directions where the electric field component is continuous.
Additionally, we assume for simplicity that the reference grid of our finite difference operator is of length one (which we can achieve by rescaling).
Further, we assume that each element has the same length $\Delta x$ in $x$- and $\Delta y$ in $y$-direction.
Both assumptions are not necessary, they are used to simplify the definition of our matrix operators.
Since we have Lagrange functions associated with boundary nodes whose support contains two elements, we define assembled quadrature weights analogous to continuous Galerkin finite element matrix assembly:
\begin{equation}
    \tilde{\omega}_i := \begin{cases}
        \omega_i \text{ for } i \in \{1,\dots,N-1\}, \\
        \omega_0 + \omega_N \text{ for } i = 0,~ i = N.
    \end{cases}
\end{equation}
In addition, we require the derivative matrix $\underline{D} \in \mathbb{R}^{N+1,N+1}$ of the finite difference operator to be of rank $N$ and $\underline{D} \underline{c} = 0$ for all constant vectors $\underline{c}$. 
Then {$\underline{V}$} defined by {Eq. \eqref{eq:vandermonde}} has full rank and the histopolation mass matrix $\underline{V}^T\underline{M}\underline{V}$ is symmetric positive definite and therefore invertible.

We can now derive our method via the Galerkin ansatz in \eqref{eq:faraday_discrete} and \eqref{eq:ampere_discrete}.
In this section we derive the scheme for the periodic case.
The non-periodic case is almost the same, but one has to adjust the index range of the sums for boundary elements and account for the fact, that some elements have no neighbour in one or more directions.
Using our basis functions as test functions we get for the first component of \eqref{eq:ampere_discrete}:
\begin{equation}
    \begin{split}
        &\sum_{k = 1}^m \sum_{i = 0}^{N-1}\sum_{j = 1}^N \int_{K_k} \bigpartialderiv{}{t}(\overline{E}^x_{i,j;k}) l^k_i(x) h^k_j(x) l^l_s(x) h^l_t(y) \text{ d}x \text{ d}y \\
        = &\sum_{k = 1}^m \sum_{i = 0}^{N-1} \sum_{j = 0}^{N-1} \int_{K_k} \bigpartialderiv{}{y}\left(B^z_{i,j;k} l^k_i(x) l^k_j(y)\right) l^l_s(x) h^l_t(y) \text{ d}x \text{ d}y\\
        = &\sum_{k = 1}^m \sum_{i = 0}^{N-1} \sum_{j = 1}^{N} \int_{K_k} \left(B^z_{i,j;k} - B^z_{i,j-1;k}\right) l^k_i(x) h^k_j(y)l^l_s(x) h^l_t(y) \text{ d}x \text{ d}y.
    \end{split}
\end{equation}
Using collocation with our finite difference quadrature rule on the grid leads us to
\begin{equation}
    \begin{split}
        &\Delta x \Delta y  \tilde{\omega}_s \sum_{j = 1}^N \bigpartialderiv{}{t}(\overline{E}^x_{s,j;l}) \sum_{o = 0}^N\omega_o {h_j^l(y^l_o)} h^l_t(y^l_o) \\
        = &\Delta x \Delta y  \tilde{\omega}_s \sum_{j = 1}^N \left(B^z_{s,j;l} - B^z_{s,j-1;l}\right) \sum_{o = 0}^N\omega_o {h_j^l(y^l_o)} h^l_t(y^l_o),
    \end{split}
\end{equation}
where the sum over the elements falls away because the support of our test functions encompasses at most two elements. We can now cancel the common factors and write it in matrix-vector form
\begin{equation}
    \underline{V}^T\underline{M}\underline{V}\frac{d}{dt}\underline{\overline{E}}_{s;l}^x = \underline{V}^T\underline{M}\underline{V} \underline{\Delta} \underline{B}_{s;l}^z,
\end{equation}
with $\frac{d}{dt}\underline{\overline{E}}_{s;l}^x := \left(\frac{d}{dt}(\overline{E}^x_{s,1;l}), \dots \frac{d}{dt}(\overline{E}^x_{s,N;l})\right)^T$ and $\underline{B}_{s;l}^z := \left(B^z_{s,0;l}, \dots B^z_{s,N;l}\right)^T$.
Usually the matrix $\underline{V}$ would contain a scaling factor $(\Delta x)^{-1}$ from the derivative matrix, but this cancels out as well.
Since the histopolation mass matrix is invertible, we obtain for all element indices $l \in \{1,\dots,m\}$
\begin{equation}
    \bigpartialderiv{}{t}(\overline{E}^x_{s,j;l}) =  B^z_{s,j;l} - B^z_{s,j-1;l},\qquad \forall s \in \{0,\dots,N\}, j \in \{1, \dots, N\}.
\end{equation}
Analogously we obtain
\begin{equation}
    \bigpartialderiv{}{t}(\overline{E}^y_{i,t;l}) =  {B^z_{i-1,t;l} - B^z_{i,t;l}},\qquad \forall t \in \{0,\dots,N\}, i \in \{1, \dots, N\}.
\end{equation}
For \eqref{eq:faraday_discrete} we obtain
\begin{equation}
    \begin{split}
        &\sum_{k = 1}^m \sum_{i = 0}^{N-1}\sum_{j = 0}^{N-1} \int_{K_k} \bigpartialderiv{}{t}({B^z_{i,j;k}}) l^k_i(x) l^k_j(x) l^l_s(x) l^l_t(y) \text{ d}x \text{ d}y \\
        = &{\sum_{k = 1}^m \sum_{i = 1}^{N} \sum_{j = 0}^{N-1} \int_{K_k} \overline{E}^y_{i,j;k} h^k_i(x) l^k_j(y) \bigpartialderiv{}{x}l^l_s(x) l^l_t(y) \text{ d}x \text{ d}y}\\
        - &{\sum_{k = 1}^m \sum_{i = 0}^{N-1} \sum_{j = 1}^{N} \int_{K_k} \overline{E}^x_{i,j;k} l^k_i(x) h^k_j(y) l^l_s(x) \bigpartialderiv{}{y}l^l_t(y) \text{ d}x \text{ d}y}.
    \end{split}
\end{equation}
Using our collocation quadrature rule and canceling scaling terms we obtain for all element indices $l \in \{1,\dots,m\}$
\begin{equation}
    \label{eq:weak_form_pointwise}
    \begin{split}
        & \tilde{\omega}_s \tilde{\omega}_t \bigpartialderiv{}{t}({B^z_{s,t;l}}) \\
        = &{\tilde{\omega}_t \sum_{i = 1}^{N} \overline{E}^y_{i,t;l} \sum_{o = 0}^N \omega_o h^l_i(x_o^l) \bigpartialderiv{}{x}l^l_s(x_o^l) + S_{s;l}}\\
        - &{\tilde{\omega}_s \sum_{j = 1}^{N} \overline{E}^x_{s,j;l} \sum_{o = 0}^N \omega_o h^l_j(y_o^l) \bigpartialderiv{}{y}l^l_t(y_o^l) + T_{t;l}},
    \end{split}
    \qquad \forall s,t \in \{0,\dots,N\}.
\end{equation}
where $S_{s;l}$ and $T_{t;l}$ denote terms from neighbouring elements if $s = 0$, $s = N$, $t = 0$, or $t = N$.
For example if $t = N$ we obtain
\begin{equation}
    T_{t;l} = \tilde{\omega}_s \sum_{j = 1}^{N} \overline{E}^x_{s,j;\sigma_u(l)} \sum_{o = 0}^N \omega_o h^{\sigma_u(l)}_j(y_o^{\sigma_u(l)}) \bigpartialderiv{}{y}l^l_t(y_o^{\sigma_u(l)}).
\end{equation}
The other terms are analogous.

\subsection{Matrix-Vector Formulation}

To assemble the matrix-vector form of the scheme, we need the 2D matrix operators. These we can obtain from the 1D element-wise operators via the Kronecker product. We repeat here the matrix operators from \ref{subsection:basis_functions}, with the special case of the derivative operators $\underline{\Delta}$ and $\underline{D}$. 
Those we have to split into two parts because the degrees of freedoms are not duplicated in the corresponding vectors.
The 1D operators also do not take into account the element size yet, this we will introduce assembling the 2D operators.
For integration we obtain the regular mass matrix $\underline{M}$ and the assembled mass matrix $\hat{\underline{M}}$:
\begin{align}
    \underline{M} = \text{diag}\left(\omega_0, \dots, \omega_N\right), \\
    \label{eq:mass_matrix_assembled}
    \underline{\hat{M}} = \text{diag}\left(\omega_0 + \omega_N, \dots, \omega_{N-1}\right).
\end{align}
The latter results from the Lagrange functions associated with boundary nodes of an element.
These have a support of two elements and are therefore integrated over both, resulting in the numerical integral $\omega_0 + \omega_N$.
We note that this case resulting in a different mass matrix is owed to the fact that we do not save duplicated degrees of freedom.
If we save all degrees of freedom of one element in one part of the overall vector and therefore include duplications in continuous directions, we could use $\underline{M}$ in all directions.
The implementation in the provided reproducibility repository does not use duplicated DOFs since we would have to properly account for the contributions of multiple elements on the duplicated DOFs due to the neighbour terms in the weak form equation for the magnetic field \eqref{eq:weak_form_pointwise}.
An element-local computation would only account for the term of the contribution of the given element, requiring inter-element communication steps.
However this is a possible option and we make no claim about which implementation is more efficient.
For our derivative operators, we have $\underline{\hat{D}} \in \mathbb{R}^{N+1 \times N}$, $\underline{\tilde{D}} \in \mathbb{R}^{N+1 \times N}$, $\underline{\hat{\Delta}} \in \mathbb{R}^{N \times N}$, and $\underline{\tilde{\Delta}} \in \mathbb{R}^{N \times N}$.
Since these are element-wise operators, the derivative at the boundaries always refers to the derivative of the function restricted to the element.
So the derivative at $x_N$ would be the left-sided derivative.
\begin{align}
    \underline{\hat{D}} &:= \left(\bigpartialderiv{}{x}l_j(x_i)\right)_{i = 0, j = 0}^{N,N-1}, \\
    \underline{\tilde{D}} &:= \begin{pmatrix}
        \bigpartialderiv{}{x}l_N(x_0) & 0 & \dots & 0 \\
        \vdots & \dots & \vdots \\
        \bigpartialderiv{}{x}l_N(x_N) & 0 & \dots & 0 \\
    \end{pmatrix}, \\
    \label{eq:curoff_delta_matrix}
    \underline{\hat{\Delta}} &:= \left(\delta_{i+1,j}-\delta_{i,j}\right)_{i = 0, j = 0}^{N-1, N-1}, \\
    \underline{\tilde{\Delta}} &:= \left(\delta_{i,N-1}\delta_{j,0}\right)_{i = 0, j = 0}^{N-1, N-1}.
\end{align}
Here $(a_{i,j})_{i = 0, j = 0}^{K,L}$ is defined as the $\mathbb{R}^{K \times L}$ matrix with entries $a_{i,j}$ with row indices running from $i = 0$ to $K$ and column indices $j = 0$ to $L$.
Since we assume $\bigpartialderiv{}{x}l_j(x_i) = D_{i,j}$ for the derivative matrix of the SBP operator from \eqref{eq:SBP_derivative_matrix} these split derivative matrices just take their value from the derivative matrix of the SBP operator.
Furthermore, we have the Vandermonde matrix $\underline{V} \in \mathbb{R}^{N+1 \times N}$ defined in \eqref{eq:vandermonde}.
To construct the $2D$ operators we need the neighbour matrices $\underline{P}^r \in \mathbb{R}^{m \times m}$ and $\underline{P}^u \in \mathbb{R}^{m \times m} $ that encode the neighbour information as follows:
\begin{align}
    P^r_{i,j} &:= \delta_{\sigma_r(i),j} \qquad \forall i,j \in \{1,\dots,m\}, \\
    P^u_{i,j} &:= \delta_{\sigma_u(i),j} \qquad \forall i,j \in \{1,\dots,m\}.
\end{align}
Let further $\underline{I}_n \in \mathbb{R}^{n \times n}$ {be} the identity matrix and $\underline{A} \otimes \underline{B}$ be the Kronecker product, then we can define the 2D mass matrices:
\begin{align}
    \underline{\hat{M}}^{2D} &:= \Delta x \Delta y \left(\underline{I_m} \otimes \underline{\hat{M}} \otimes \underline{\hat{M}}\right), \\
    \underline{M}_x^{2D} &:= \Delta x \Delta y \left(\underline{I_m} \otimes \underline{M} \otimes \underline{\hat{M}}\right), \\
    \underline{M}_y^{2D} &:= \Delta x \Delta y \left(\underline{I_m} \otimes \underline{\hat{M}} \otimes \underline{M}\right), \\
    \underline{M}^{2D} &:= \Delta x \Delta y \left(\underline{I_m} \otimes \underline{M} \otimes \underline{M}\right).
\end{align}
We also obtain the 2D derivative operators:
\begin{align}
    \underline{D}^{2D}_x &:= \left(\Delta x\right)^{-1} \left(\underline{I}_{m} \otimes \underline{I}_N \otimes \hat{\underline{D}} + \underline{P}_r \otimes \underline{I}_N \otimes \tilde{\underline{D}}\right), \\
    \underline{D}^{2D}_y &:= \left(\Delta y\right)^{-1} \left(\underline{I}_{m} \otimes \hat{\underline{D}} \otimes \underline{I}_N + \underline{P}_u \otimes \tilde{\underline{D}} \otimes \underline{I}_N\right), \\
    \underline{\Delta}^{2D}_x &:= \underline{I}_{m} \otimes \underline{I}_N \otimes \hat{\underline{\Delta}} + \underline{P}_r \otimes \underline{I}_N \otimes \tilde{\underline{\Delta}}, \\
    \underline{\Delta}^{2D}_y &:= \underline{I}_{m} \otimes \hat{\underline{\Delta}} \otimes \underline{I}_N + \underline{P}_u \otimes \tilde{\underline{\Delta}} \otimes \underline{I}_N.
\end{align}
Since these operators applied to the variable vectors yield the exact derivative, we have again
\begin{equation}
    \label{eq:discrete_schwarz}
    \underline{\Delta}^{2D}_x \underline{\Delta}^{2D}_y - \underline{\Delta}^{2D}_x \underline{\Delta}^{2D}_y = 0.
\end{equation}
The Vandermonde matrices are
\begin{align}
    \underline{V}^{2D}_x &:= \left(\Delta x\right)^{-1}\left(\underline{I}_{m} \otimes \underline{I}_{N} \otimes \underline{V}\right),\\
    \underline{V}^{2D}_y &:= \left(\Delta y\right)^{-1}\left(\underline{I}_{m} \otimes \underline{V} \otimes \underline{I}_{N}\right).
\end{align}
We note that we again have the analogous equations to \eqref{eq:derivative_vandermonde_equality}
\begin{align}
    \label{eq:derivative_vandermonde_equality_x}
    \underline{D}_x &= \underline{V}^{2D}_x \underline{\Delta}^{2D}_x, \\
    \label{eq:derivative_vandermonde_equality_y}
    \underline{D}_y &= \underline{V}^{2D}_y \underline{\Delta}^{2D}_y,    
\end{align}
which leads to the semi-discrete scheme
\begin{align}
    \label{eq:semidiscrete_ampere_x}
    \frac{\text{d}}{\text{d}t} \overline{\underline{E}}^x &= \underline{\Delta}^{2D}_y \underline{B}^z, \\
    \label{eq:semidiscrete_ampere_y}
    \frac{\text{d}}{\text{d}t} \overline{\underline{E}}^y &= -\underline{\Delta}^{2D}_x \underline{B}^z,\\ 
    \label{eq:semidiscrete_faraday}
    \underline{\hat{M}}^{2D} \frac{\text{d}}{\text{d}t} \underline{B}^z &= \underline{\Delta}^{2D,T}_x \underline{V}_x^{2D,T} \underline{M}_y^{2D} \underline{V}^{2D}_x \underline{\overline{E}}^y - \underline{\Delta}^{2D,T}_y \underline{V}^{2D,T}_y \underline{M}_x^{2D} \underline{V}^{2D}_y \underline{\overline{E}}^x.
\end{align}
If we want to impose homogeneous boundary conditions, we have to change the matrix operators slightly. 
We need versions for the assembled mass matrix \eqref{eq:mass_matrix_assembled} and the difference matrix \eqref{eq:delta_matrix} for the left/lower and the right/upper boundary elements.
The assembled mass matrix for the boundary elements read
\begin{align}
    \underline{\hat{M}}_l &= \text{diag}\left(\omega_0, \dots, \omega_{N-1}\right),\\
    \underline{\hat{M}}_r &= \text{diag}\left(\omega_0 + \omega_N, \dots, \omega_{N-1}, \omega_N\right),\\
    \underline{\hat{M}}_1 &= \text{diag}\left(\omega_0, \dots, \omega_{N-1}, \omega_N\right),
\end{align}
where $\underline{\hat{M}}_1$ is the assembled mass matrix in the case of one element.
The boundary replacements for the cutoff difference matrices $\hat{\underline{\Delta}}$ read
\begin{align}
    \underline{\hat{\Delta}}_l &= \underline{\hat{\Delta}}, \\
    \underline{\hat{\Delta}}_r &= \underline{\Delta}, \\
    \underline{\hat{\Delta}}_1 &= \underline{\Delta}.
\end{align}
The remaining part $\underline{\tilde{\Delta}}$ of the difference matrix entries remain the same, but the dimensions have to be adjusted to match the number of DOFs in the boundary elements by adding additional zeros.

The 2D matrices do no longer have the form of a Kronecker product over all elements.
The mass and Vandermonde matrices are still block-diagonal matrices with one block for each element, but since we have to adjust for the boundary elements, we cannot take a Kronecker product with $I_m$ anymore.
For the 2D Vandermonde matrix the only thing that changes for different blocks is the size of the identity matrix in the Kronecker product.
For example, if we have an element $K$ on the upper boundary, then the block for $V_x$ for this element reads
\begin{equation}
    V^{2D}_{x,K} = \underline{I}_{N+1} \otimes \underline{V}.
\end{equation}
For the 2D mass and Vandermonde matrices we have to use the boundary versions for $\hat{M}$, $\underline{\hat{\Delta}}$ and $\underline{\tilde{\Delta}}$.
The mass matrix $\underline{M}^{2D}$, that only uses element-local mass-matrices, remains unchanged.
If we want to impose essential boundary conditions, we have to truncate the boundary matrices, since the test functions corresponding to boundary DOFs are no longer part of the test space.
Then we have
\begin{align}
    \underline{\hat{M}}_l &= \text{diag}\left(\omega_1, \dots, \omega_{N-1}\right),\\
    \underline{\hat{M}}_r &= \text{diag}\left(\omega_0 + \omega_N, \dots, \omega_{N-1}\right),\\
    \underline{\hat{M}}_1 &= {\text{diag}\left(\omega_1, \dots, \omega_{N-1} \right)},
\end{align}
and
\begin{align}
    \underline{\hat{\Delta}}_l &= {\underline{\Delta}_{1:N,2:N}}, \\
    \underline{\hat{\Delta}}_r &= \underline{\Delta}_{1:N,1:N}, \\
    \underline{\hat{\Delta}}_1 &= \underline{\Delta}_{1:N,2:N},
\end{align}
where $\underline{\Delta}_{1:N, 1:N}$ is the matrix we get when we take all rows, but only the first $N$ columns from $\underline{\Delta}$.
Alternatively one can use the untruncated versions and set the boundary DOFs to $0$ after each right-hand side evaluation, which is the variant that is implemented in the linked github repository.
Using the Summation-by-Parts property we can also derive a strong form version of the scheme.
First we rewrite the semi-discrete Faraday equation \eqref{eq:semidiscrete_faraday}, using the identity $\underline{V}\underline{\Delta} = \underline{D}$
\begin{equation}
    \label{eq:semidiscrete_faraday_kronecker}
    \begin{split}
        \frac{\text{d}}{\text{d}t} \underline{B}^z = & \frac{1}{\Delta x^2} \left(\underline{I}_m \otimes \underline{I}_N \otimes \underline{\hat{M}}^{-1} \underline{\hat{D}}^T \underline{M} \underline{V} 
        + \underline{P}_r^T \otimes \underline{I}_N \otimes \underline{\hat{M}}^{-1} \underline{\tilde{D}}^T \underline{M} \underline{V}\right) \underline{\overline{E}}^y \\
        - & \frac{1}{\Delta y^2}\left(\underline{I}_m \otimes \underline{\hat{M}}^{-1} \underline{\hat{D}}^T \underline{M} \underline{V} \otimes \underline{I}_N
        +\underline{P}_u^T \otimes  \underline{\hat{M}}^{-1} \underline{\tilde{D}}^T \underline{M} \underline{V} \otimes \underline{I}_N\right) \underline{\overline{E}}^x.
    \end{split}
\end{equation}
To deal with the split derivative matrix, we will denote for a $N-1 \times N-1$ matrix $\underline{A}$ by $\underline{A}_{0:N-1} \in M^{N+1 \times N+1}(\mathbb{R})$ the matrix with the last row cut off and by $\underline{A}_{N;0} \in M^{N+1 \times N+1}(\mathbb{R})$ the matrix having as first row the last row of $\underline{A}$, with all other entries being $0$.
Applying the Summation-by-Parts formula \eqref{eq:SBP} now, we obtain for the first term of the {right} side of \eqref{eq:semidiscrete_faraday_kronecker}:

\begin{equation}
    \begin{split}
        & \frac{1}{\Delta x^2}\left(\underline{I}_m \otimes \underline{I}_N \otimes \underline{\hat{M}}^{-1} \underline{\hat{D}}^T \underline{M} \underline{V}
        + \underline{P}_r^T \otimes \underline{I}_N \otimes \underline{\hat{M}}^{-1} \underline{\tilde{D}}^T \underline{M} \underline{V}\right)\underline{\overline{E}}^y \\
        = &\frac{1}{\Delta x^2}\left(\underline{I}_m \otimes \underline{I}_N \otimes \underline{\hat{M}}^{-1} \underline{B}_{0:N-1}\underline{V}
        + \underline{P}_r^T \otimes \underline{I}_N \otimes \underline{\hat{M}}^{-1} \underline{B}_{N;0}  \underline{V}\right)\underline{\overline{E}}^y \\
        - &\frac{1}{\Delta x^2} \left(\underline{I}_m \otimes \underline{I}_N \otimes \underline{\hat{M}}^{-1} \underline{M}_{0:N-1} \underline{D} \underline{V}
        +\underline{P}_r^T \otimes \underline{I}_N \otimes \underline{\hat{M}}^{-1} \underline{M}_{N;0} \underline{D} \underline{V}\right)\underline{\overline{E}}^y.
    \end{split}
\end{equation}

An analogous equation holds for the second term.
Component wise we obtain
\begin{equation}
    \label{eq:strong_form_componentwise}
    \begin{split}
        &\frac{\text{d}}{\text{d}t} B^z_{k,l,m} = - (1 - \delta_{0,k}) \frac{1}{\Delta x}\sum_{i = 0}^N D_{k,i} E^y_{i,l,m} + \frac{1}{\Delta y}(1 - \delta_{0,l}) \sum_{i = 0}^N D_{l,i} E^x_{k,i,m} \\
        &-\delta_{0,k} \frac{1}{\Delta x}\frac{1}{\omega_N + \omega_0} \left(\omega_N \sum_{i = 0}^N D_{N,i} E^y_{i,l,\sigma_{le}(m)} + \omega_0 \sum_{i = 0}^N D_{1,i} E^y_{i,l,m} \right) \\
        &+ \delta_{0,l} \frac{1}{\Delta y}\frac{1}{\omega_N + \omega_0} \left(\omega_N \sum_{i = 0}^N D_{N,i} E^x_{k,i,\sigma_{lo}(m)} + \omega_0 \sum_{i = 0}^N D_{1,i} E^x_{k,i,m} \right) \\
        &+ \delta_{0,k}  \frac{1}{\Delta x} \frac{1}{\omega_N + \omega_0} \left(E^y_{N,l,\sigma_{le}(m)} - E^y_{0,l,m}\right) \\
        &- \delta_{0,l} \frac{1}{\Delta y} \frac{1}{\omega_N + \omega_0} \left(E^x_{k,N,\sigma_{lo}(m)} - E^x_{k,0,m} \right),
    \end{split}
\end{equation}
where $E^x_{k,l,m}$ and $E^x_{k,l,m}$ are the nodal values of the electric field computed with the scaled Vandermonde matrix while $\sigma_{lo}(m)$ and $\sigma_{le}(m)$ are the indices of the element below and left of the element with index $m$. 
Since we already use the scaled Vandermonde matrix in the nodal evaluation, the scaling factors in the component-wise equation are not squared anymore.

We note that on a function level, the derivative matrix $\underline{D}$ computes the nodal derivatives of the element-local Lagrange functions, not the histopolation functions.
Since we take the derivative in the directions we have the histopolation basis functions, this amounts to interpolating the functions with the element-local Lagrange basis, without requiring continuity, and then computing the local derivative.
The strong form equation \eqref{eq:strong_form_componentwise} then consists of the element-local interpolation derivative in x-direction for $E^y$ and in y-direction for $E^x$ for $k \neq 0$ and $l \neq 0$ respectively. 
If the degree of freedom in question is on the element boundary in {the} x- or y-direction, we instead have a weighted jump penalty term minus the weighted average of the discontinuous derivative at the boundary.
This results in a treatment of the discontinuous derivative similar to a finite difference SAT, a finite volume or a discontinuous Galerkin ansatz.

Proceeding from the weak formulation, we can rewrite the scheme using the discrete Hamiltonian or discrete energy defined as
\begin{equation}
    \label{eq:hamiltonian}
    \begin{split}
        H(\underline{U}) = \frac{1}{2} (&\overline{\underline{E}}^{x,T} \underline{V}^{2D,T}_y M_x^{2D} \underline{V}^{2D}_y \overline{\underline{E}}^x \\
        + &\overline{\underline{E}}^{y,T} \underline{V}^{2D,T}_x \underline{M}_y^{2D} \underline{V}^{2D}_x \overline{\underline{E}}^y + \underline{B}^{z,T} \underline{\hat{M}}^{2D} \underline{B}^z),
    \end{split}
\end{equation}
where
\begin{equation}
    \underline{U} = \begin{pmatrix}
        \overline{\underline{E}}^{x} \\
        \overline{\underline{E}}^{y} \\
        \underline{B}^z
    \end{pmatrix},
\end{equation}
and
\begin{equation}
    \label{eq:poisson}
    \underline{J} = \begin{pmatrix}
        0 & 0 & \underline{\Delta}_y^{2D} \underline{\hat{M}}^{2D,-1} \\
        0 & 0 & -\underline{\Delta}_x^{2D} \underline{\hat{M}}^{2D,-1} \\
        -\underline{\hat{M}}^{2D,-1} \underline{\Delta}_y^{2D,T} & \underline{\hat{M}}^{2D,-1} \underline{\Delta}_x^{2D,T} & 0
    \end{pmatrix},
\end{equation}
is the Poisson matrix.

Then the scheme \eqref{eq:semidiscrete_ampere_x}, \eqref{eq:semidiscrete_ampere_y}, and \eqref{eq:semidiscrete_faraday} is equivalent to
\begin{equation}
    \label{eq:hamiltonian_form}
    \frac{\text{d}}{\text{d}t} \underline{U} = \underline{J} \nabla_U H(\underline{U}),
\end{equation}
with
\begin{equation}
    \nabla_U H(\underline{U}) = \begin{pmatrix}
        \underline{V}^{2D,T}_y \underline{M}_x^{2D} \underline{V}^{2D}_y \overline{\underline{E}}^{x} \\
        \underline{V}^{2D,T}_x \underline{M}_y^{2D} \underline{V}^{2D}_x \overline{\underline{E}}^{y} \\
        \underline{\hat{M}}^{2D} \underline{B}^z
    \end{pmatrix}.
\end{equation}
The semi-discrete energy conservation of the scheme follows directly:
\begin{equation}
    \frac{\text{d}}{\text{d}t} H(\underline{U}) = \left(\nabla_U H(\underline{U})\right)^T \frac{\text{d}}{\text{d}t}\underline{U} =  \left(\nabla_U H(\underline{U})\right)^T \underline{J} \left(\nabla_U H(\underline{U})\right) = 0,
\end{equation}
since $\underline{J}$ is anti-symmetric.
We remark that this property can only be proven for the weak form of the scheme.
For non-SBP operators {the} strong and weak form are not equivalent, so this may not hold for the strong form of the scheme for non-SBP operators.

\subsection{Time Integration}

We will consider two different time integration schemes here.
One is an explicit three-stage Runge-Kutta scheme \cite{shu1988} and the other, the discrete gradient average vector field method \cite{mclachan1999}, which in this case coincides with the Crank-Nicolson scheme \cite{crank_nicolson_1947}, which is implicit.
The purpose of the latter is to exactly preserve the discrete energy in the fully discrete setting.
The ansatz for the discrete gradient method is to find an approximation $\overline{\nabla}_U H(\underline{U}_{n}, \underline{U}_{n+1})$, that depends on the current and the next time step, for $\nabla H(\underline{U})$.
We require this approximation to fulfill
\begin{equation}
    \label{eq:discrete_gradient_condition}
    \left(\underline{U}_{n+1} - \underline{U}_{n} \right)^T \overline{\nabla}_U H(\underline{U}_{n}, \underline{U}_{n+1}) = H(\underline{U}_{n+1}) - H(\underline{U}_{n}),
\end{equation}
because definining the time integration scheme implicitly by
\begin{equation}
    \frac{1}{\Delta t} \left(\underline{U}_{n+1} - \underline{U}_{n} \right) = \underline{J} \overline{\nabla}_U H(\underline{U}_{n}, \underline{U}_{n+1}),
\end{equation}
conserves the Hamiltonian \cite{mclachan1999} since
\begin{equation}
    \begin{split}
        &H(\underline{U}_{n+1}) - H(\underline{U}_{n}) = \left(\underline{U}_{n+1} - \underline{U}_{n} \right)^T \overline{\nabla}_U H(\underline{U}_{n}, \underline{U}_{n+1}) \\
        = &\Delta t \left(\overline{\nabla}_U H\right)^T \underline{J}^T \left(\overline{\nabla}_U H\right) = 0.
    \end{split}
\end{equation}
In this case the average vector field discrete gradient is, 
\begin{equation}
    \begin{split}
    \overline{\nabla}_UH(\underline{U}_n,\underline{U_{n+1}}) &\coloneqq \int_0^1 \left(\lambda  \nabla_U H(\underline{U}_{n}) + (1 -\lambda) \nabla_U H(\underline{U}_{n+1})\right) \text{d}\lambda \\
    &= \frac{1}{2} \nabla_U H(\underline{U}_{n+1} + \underline{U}_n).
    \end{split}
\end{equation}
Since $\nabla_UH$ is linear, this choice leads to the Crank-Nicolson time integration method \cite{crank_nicolson_1947}:

\begin{equation}
    \frac{1}{\Delta t} \left( 
    \underline{U}_{n+1}
    -
    \underline{U}_{n}
    \right)
    = \frac{1}{2} \underline{J} \nabla_U H(\underline{U}_{n+1} + \underline{U}_n).
\end{equation}

\begin{proposition}
    Both the Crank-Nicolson method and any Runge-Kutta scheme applied to the semi-discrete method given by \eqref{eq:semidiscrete_ampere_x}, \eqref{eq:semidiscrete_ampere_y} and \eqref{eq:semidiscrete_faraday} exactly preserve the divergence of the electric field in time. If the divergence was $0$ at the start time $t_0$, it remains $0$ at all time-steps $t_n$.
\end{proposition}
\begin{proof}
    The update step for a $s$-stage Runge-Kutta scheme for the electric field has the form
    \begin{equation}
        \begin{pmatrix}
            \overline{\underline{E}}_{n+1}^{x} \\
            \overline{\underline{E}}_{n+1}^{y}
        \end{pmatrix}
        =
        \begin{pmatrix}
            \overline{\underline{E}}_{n}^{x} \\
            \overline{\underline{E}}_{n}^{y}
        \end{pmatrix}
        + \Delta t \begin{pmatrix}
            \underline{\Delta}_y^{2D} \\
            -\underline{\Delta}_x^{2D}
        \end{pmatrix}
        \sum_{j = 1}^s {\underline{k}_j},
    \end{equation}
    for the vectors {$\underline{k}_j$} that are generally dependent on all other stages of the Runge-Kutta scheme.
    Applying the discrete divergence on both sides yields
    \begin{equation}
        \begin{split}
            &\begin{pmatrix}
                \underline{\Delta}_x^{2D} \\
                \underline{\Delta}_y^{2D}
            \end{pmatrix}
            \cdot
            \begin{pmatrix}
                \overline{\underline{E}}_{n+1}^{x} \\
                \overline{\underline{E}}_{n+1}^{y}
            \end{pmatrix}
            =
            \begin{pmatrix}
                \underline{\Delta}_x^{2D} \\
                \underline{\Delta}_y^{2D}
            \end{pmatrix}
            \cdot
            \begin{pmatrix}
                \overline{\underline{E}}_{n}^{x} \\
                \overline{\underline{E}}_{n}^{y}
            \end{pmatrix}
            \\
            + &\Delta t 
                \begin{pmatrix}
                \underline{\Delta}_x^{2D} \\
                \underline{\Delta}_y^{2D}
            \end{pmatrix}
            \cdot
            \begin{pmatrix}
                \underline{\Delta}_y^{2D} \\
                -\underline{\Delta}_x^{2D}
            \end{pmatrix}
            \sum_{j = 1}^s \underline{k}_j
            = \begin{pmatrix}
                \underline{\Delta}_x^{2D} \\
                \underline{\Delta}_y^{2D}
            \end{pmatrix}
            \cdot
            \begin{pmatrix}
                \overline{\underline{E}}_{n}^{x} \\
                \overline{\underline{E}}_{n}^{y}
            \end{pmatrix}
        \end{split}
    \end{equation}
    where we used \eqref{eq:discrete_schwarz} and
    \begin{equation}
            \begin{pmatrix}
                \underline{A} \\
                \underline{B}
            \end{pmatrix}
            \cdot
            \begin{pmatrix}
                \underline{v} \\
                \underline{w}
            \end{pmatrix}
            \vcentcolon =
            \begin{pmatrix}
                \underline{A}\underline{v} \\
                \underline{B}\underline{w}
            \end{pmatrix},
    \end{equation}
    for matrices $\underline{A}, \underline{B}$ and vectors $\underline{v}, \underline{w}$ with compatible dimensions.
    For the Crank-Nicolson method we have
    \begin{equation}
        \label{eq:mixed_time_step_formulation}
        \frac{1}{\Delta t} \left( 
            \begin{pmatrix}
            \underline{\overline{E}}_{n+1}^x\\
            \underline{\overline{E}}_{n+1}^y\\
            \underline{B}^z_{n+1}
        \end{pmatrix}
        -
         \begin{pmatrix}
            \underline{\overline{E}}_n^x\\
            \underline{\overline{E}}_n^y\\
            \underline{B}^z_n
        \end{pmatrix}
        \right)
        = \frac{1}{2} \underline{J} \nabla_U H \left(
            \begin{pmatrix}
                \underline{\overline{E}}_{n+1}^x\\
                \underline{\overline{E}}_{n+1}^y\\
                \underline{B}^z_{n+1}
            \end{pmatrix}
            +
            \begin{pmatrix}
                \underline{\overline{E}}_n^x\\
                \underline{\overline{E}}_n^y\\
                \underline{B}^z_n
            \end{pmatrix}        
        \right).
    \end{equation}
    Applying the discrete divergence to the electric field yields
    \begin{equation}
        \frac{1}{\Delta t}
        \begin{pmatrix}
            \underline{\Delta}_x^{2D} \\
            \underline{\Delta}_y^{2D} \\
            0
        \end{pmatrix}
        \cdot
        \left( 
        \begin{pmatrix}
            \underline{\overline{E}}_{n+1}^x\\
            \underline{\overline{E}}_{n+1}^y\\
            \underline{B}_{n+1}^z
        \end{pmatrix}
        -
        \begin{pmatrix}
            \underline{\overline{E}}_n^x\\
            \underline{\overline{E}}_n^y\\
            \underline{B}_n^z
        \end{pmatrix}
        \right)
        =
        0,
    \end{equation}
    or
    \begin{equation}
        \begin{pmatrix}
            \underline{\Delta}_x^{2D} \\
            \underline{\Delta}_y^{2D} \\
            0
        \end{pmatrix}
        \cdot
        \begin{pmatrix}
            \underline{\overline{E}}_{n+1}^x\\
            \underline{\overline{E}}_{n+1}^y\\
            \underline{B}_{n+1}^z
        \end{pmatrix}
        =
        \begin{pmatrix}
            \underline{\Delta}_x^{2D} \\
            \underline{\Delta}_y^{2D} \\
            0
        \end{pmatrix}
        \cdot
        \begin{pmatrix}
            \underline{\overline{E}}_n^x\\
            \underline{\overline{E}}_n^y\\
            \underline{B}_n^z
        \end{pmatrix},
    \end{equation}
    which completes the proof.
\end{proof}

\section{Numerical Verification}

We look at two test cases.
The first is from \cite{sonnendruecker2012} with domain $\Omega = [-1,1]^2$, end time $T = 1$, periodic boundary conditions, and the initial condition
\begin{align}
    E^x(x,y) &= 0,\\
    E^y(x,y) &= 0,\\
    B^z(x,y) &= \text{cos}(\pi x + \pi) \text{cos}(\pi y + \pi),
\end{align}
with the analytical solution
\begin{align}
    E^x(x,y,t) &= -\frac{1}{\sqrt{2}} \text{cos}(\pi x + \pi) \text{sin}(\pi y + \pi) \text{sin}(\sqrt{2} \pi t),\\
    E^y(x,y,t) &= \frac{1}{\sqrt{2}} \text{sin}(\pi x + \pi) \text{cos}(\pi y + \pi) \text{sin}(\sqrt{2} \pi t),\\
    B^z(x,y,t) &= \text{cos}(\pi x + \pi) \text{cos}(\pi y + \pi) \text{cos}(\sqrt{2} \pi t).
\end{align}
The second test case has domain $\Omega = [0,1]^2$, end time $T = 1$, initial condition
\begin{align}
    E^x(x,y) &= 0,\\
    E^y(x,y) &= 0,\\
    B^z(x,y) &= \sqrt{2} \text{cos}(\pi x) \text{sin}(\pi y),
\end{align}
with essential homogeneous boundary condition for $y = 0$ and $y = 1$, and natural boundary conditions for $x = 0$ and $x = 1$, resulting in:
\begin{align}
    E^y(x,y,t) &= 0 \text{ on } \partial \Omega \times [0, 1],\\
    {B^z(x,0,t)} &\equiv 0, \\
    {B^z(x,1,t)} &\equiv 0.
\end{align}
This problem has the analytical solution
\begin{align}
    E^x(x,y,t) &= \text{cos}(\pi x) \text{cos}(\pi y) \text{sin}(\sqrt{2} \pi t),\\
    E^y(x,y,t) &= \text{sin}(\pi x) \text{sin}(\pi y) \text{sin}(\sqrt{2} \pi t),\\
    B^z(x,y,t) &= \sqrt{2} \text{cos}(\pi x) \text{sin}(\pi y) \text{cos}(\sqrt{2} \pi t).
\end{align}
For our tests we use a Julia code, using scripts from \cite{fernandez2014review}. To solve the linear system for the Crank-Nicolson method, we use the GMRES solver \cite{saad1986} implemented in the Julia package Krylov.jl \cite{montoison-orban-2023}. For the repository, see the data availability section. The tests are run with Julia 1.11.7 on {Ubuntu 24.04.4 LTS} with a AMD Ryzen 9 9900X processor.

\subsection{Convergence Tests}

We run two different kinds of convergence tests, where we use the discrete $L^2$ norm induced by our global mass matrices. 
This means we compute the error between the nodal values of the numerical approximation $\underline{E}^x_{\text{num}}$, $\underline{E}^y_{\text{num}}$, $\underline{B}^z_{\text{num}}$ and the exact nodal values at time $T = 1$, $\underline{E}^x$, $\underline{E}^y$, $\underline{B}^z$ as follows:
\begin{align}
    \|\underline{E}^x_{\text{num}} - \underline{E}^x\|_{L^2, E_x, num} &\vcentcolon = \sqrt{\left(\underline{E}^x_{\text{num}} - \underline{E}^x\right)^T \underline{M}^{2D}_x \left(\underline{E}^x_{\text{num}} - \underline{E}^x\right)},\\
    \|\underline{E}^y_{\text{num}} - \underline{E}^y\|_{L^2, E_y, num} &\vcentcolon = \sqrt{\left(\underline{E}^y_{\text{num}} - \underline{E}^y\right)^T \underline{M}^{2D}_y \left(\underline{E}^y_{\text{num}} - \underline{E}^y\right)},\\
    \|\underline{B}^z_{\text{num}} - \underline{B}^z\|_{L^2, B^z, num} &\vcentcolon = \sqrt{\left(\underline{B}^z_{\text{num}} - \underline{B}^z\right)^T \underline{\hat{M}}^{2D} \left(\underline{B}^z_{\text{num}} - \underline{B}^z\right)}.
\end{align}
For the first kind of convergence test we keep the number of points per element constant but vary the amount of elements.
For time integration we use a three-stage strong stability-preserving Runge-Kutta scheme (see \cite{shu1988}, equation (2.18)) with a constant time step of $\Delta t = 2 \cdot 10^{-5}$. 
We run this test for the operator of boundary order $3$ by Strand \cite{strand1994} for 12 and 13 points per element in the periodic case and for 12 points again for the non-periodic case to investigate the behavior of operators dominated by boundary stencils.
Appendix B contains further results the periodic case for 14 points per element and similar results for the operator of boundary order $2$ by Strand \cite{strand1994}.
Because of the symmetry of the example, the errors for $E^x$ and $E^y$ are identical up to machine precision and are therefore only listed once in the periodic case.

We observe an even-odd behavior in Tables \ref{tbl:eoc_p3_12n} and \ref{tbl:eoc_p3_13n} and further results in appendix B where we obtain an alternating order of $p$ and $p+1$ for the electric field, depending on the number of points per element.
This behavior is reminiscent of the even-odd behavior previously reported for central (non-dissipative) continuous and discontinuous Galerkin discretizations, see for instance \cite{hindenlang2020} and \cite{angulo2015nodal}.
However, the analogy should be interpreted carefully. 
In the Galerkin methods, the discretization is parameterized by the polynomial degree, whereas the present formulation is defined directly through nodal SBP finite-difference operators and therefore has no underlying polynomial-degree parameter. 
Consequently, the analogous behavior appears naturally in the number of nodal degrees of freedom per element.

As detailed above, present method employs a continuous-Galerkin-type superposition at element interfaces rather than dissipative SAT penalties or upwind numerical fluxes. 
The resulting discretization is therefore essentially non-dissipative, similarly to central continuous and discontinuous Galerkin formulations where parity-dependent convergence behavior has previously been observed \cite{hindenlang2020,angulo2015nodal}.
While we are not aware of previous studies reporting the exact same behavior for central SBP finite-difference discretizations, related even-odd effects have been reported for SBP finite-difference operators in connection with the treatment of high-frequency numerical modes \cite{mattsson2017diagonal,duru2020upwind}. This suggests that the observed behavior may be linked to the discrete high-frequency mode structure of the nodal SBP operators and the absence of interface dissipation, rather than to the formal truncation order alone.

\begin{figure}
    \centering
    \includegraphics[width=0.8\linewidth]{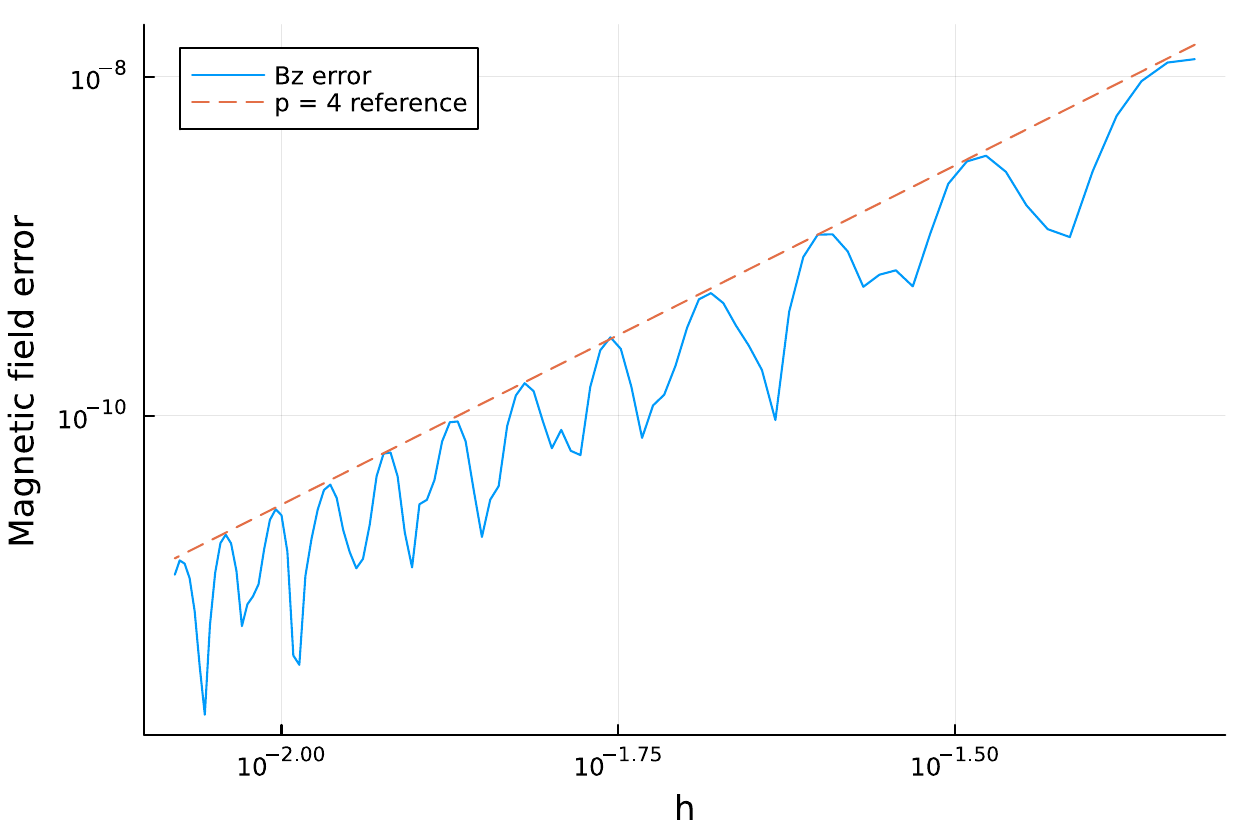}
    \caption{Log-log plot of the magnetic field error against the element size for the non-periodic problem for degree $p = 3$ with SSPRK time integration with element number between $21$ and $120$ elements}
    \label{fig:convergence_study_B_np}
\end{figure}

The magnetic field EOC varies without clearly converging towards a set number, though the error of the magnetic field tends to improve at least as good if not faster then the electric field error.
The non-periodic case for $12$ points in table \ref{tbl:eoc_p3_12n_np} shows the same behavior as the periodic case.
To further investigate the behavior of the magnetic field convergence, we ran additional simulations for the non-periodic case with the $p = 3$ operator with a varying amount of elements.
We ran the problem to end time $T = 0.1$ with element numbers per direction from $21$ to $120$ elements.
Figure \ref{fig:convergence_study_B_np} shows a log-log plot of element width $h$ against the discrete {$L^2$}-error of the magnetic field.
We see that the convergence of the magnetic field is bounded from above by a polynomial of degree $4$, but has a highly non-monotonous convergence behavior, exhibiting a periodic pattern (when seen as a function in $h^{-1}$).
Further analysis of this behavior is beyond the scope of this paper, but we see that at least in this case we have a convergence behavior that is still $\mathcal{O}(h^4)$.
The second kind of test is again run using operators from \cite{strand1994} of polynomial degree $3$ for both the periodic and the non-periodic case. 
We again use the SSPRK time integration scheme, but for the periodic case, further results using the Crank-Nicolson time integration scheme and also with the $p = 2$ operator can be found in appendix B and exhibit similar behavior.
We again use a constant time step of $\Delta t = 2 \cdot 10^{-5}$. 
For the spatial discretization we just use one element, but a varying number of nodes per element.
In tables \ref{tbl:eoc_p3_ssprk} and \ref{tbl:eoc_p3_ssprk_np} we see higher polynomial orders than in the first test, which is likely due to the fact that the inner stencils, which are of degree $2p$ become more influential and increase the EOC.

We also ran this convergence test for the periodic test case with a non-SBP operator using the mass matrix of the $p = 2$ and the derivative matrix of the $p = 3$ operator, which together do not have the SBP property.
Since {the strong and weak forms} of the scheme are no longer equivalent for this operator, we run a convergence test for both the strong and weak forms of the scheme.
Table \ref{tbl:eoc_non_sbp_strong_ssprk} shows the convergence rates for the strong from, exhibiting similar convergence rates to the $p = 3$ SBP operator test case, while the convergence rates of the weak form variant in table \ref{tbl:eoc_non_sbp_weak_ssprk} are significantly lower.
\subsection{Divergence Tests}

\begin{figure}[H]
    \centering
    \includegraphics[width=0.8\linewidth]{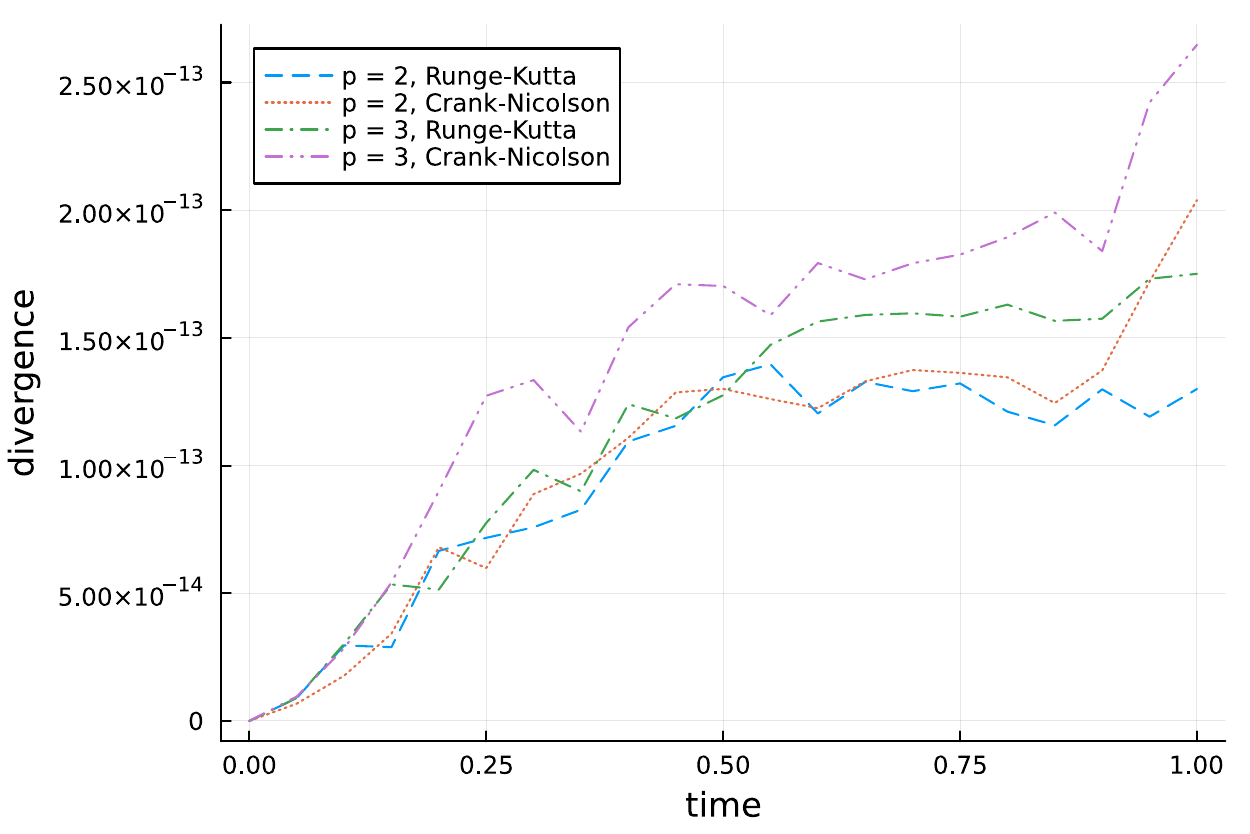}
    \caption{Divergence errors for the periodic problem for operators of degree $p = 2,3$ with both the SSPRK and Crank-Nicolson time integration, $5$ elements with $20$ points per dimension}
    \label{fig:plot_div_fine_grid}
\end{figure}

\begin{figure}[H]
    \centering
    \includegraphics[width=0.8\linewidth]{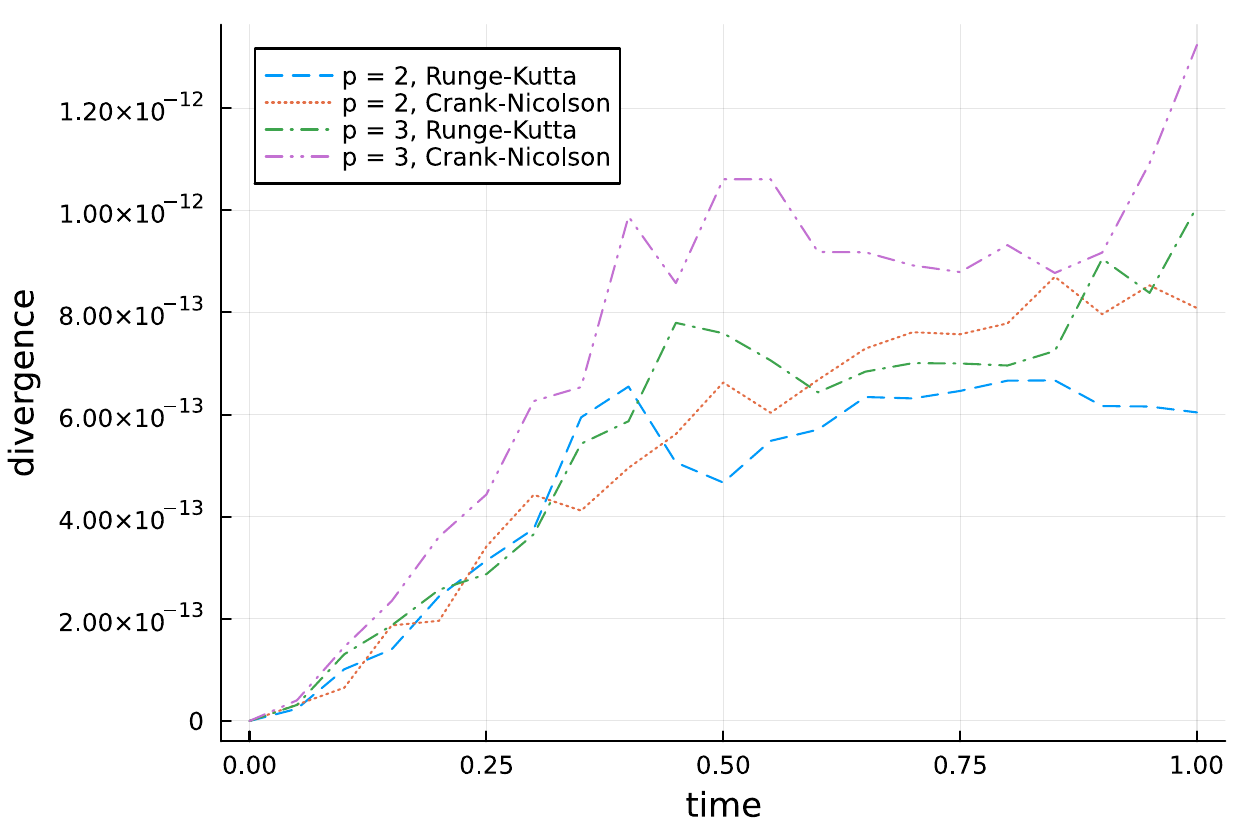}
    \caption{Divergence errors for the non-periodic problem for operators of degree $p = 2,3$ with both the SSPRK and Crank-Nicolson time integration, $5$ elements with $20$ points per dimension}
    \label{fig:plot_div_fine_grid_np}
\end{figure}

A key property of the proposed discretization is the exact preservation of the divergence condition \eqref{eq:gauss_law}. 
We test this using both a discretization using $5$ elements with $20$ points in each direction for both the periodic and the non-periodic test case.
We test divergence preservation with the operators of degree $2$ and $3$ and both with the SSPRK and the Crank-Nicolson scheme.
The time step is defined by
\begin{equation}
    \label{eq:cfl}
    \Delta t = \text{CFL} \hspace{0.1cm} \omega_{min},
\end{equation}
for $\text{CFL} = 1$. 
Here $\omega_{min}$ is the minimal quadrature weight of the global mass matrix $\hat{M}_{2D}$.
Because we assume unit light speed, the maximum absolute value of the wave speeds is $1$ and therefore does not appear in \eqref{eq:cfl}.
The figures \ref{fig:plot_div_fine_grid} and \ref{fig:plot_div_fine_grid_np} show the maximum absolute value of the divergence of the electric field at the grid points.
We can observe that the solution remains divergence free up to machine precision during the entire simulation, regardless of operator, time integration scheme, or grid size.
This confirms the previously proven fully discrete divergence preservation.

\subsection{Energy Conservation Tests}

\begin{figure}[H]
    \centering
    \includegraphics[width=0.8\linewidth]{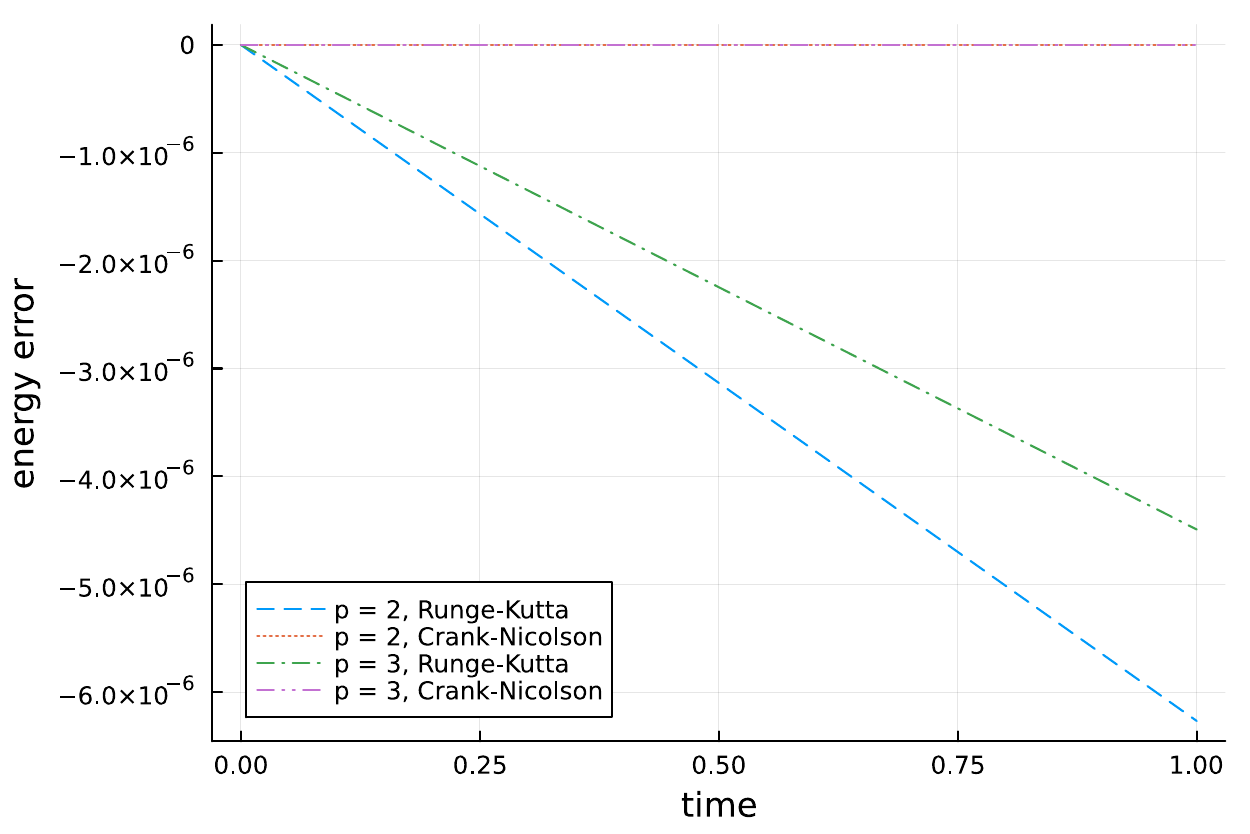}
    \caption{Energy errors for the periodic problem for operators of degree $p = 2,3$ with both the SSPRK and Crank-Nicolson time integration, $5$ elements with $20$ points per dimension}
    \label{fig:plot_energy_fine_grid}
\end{figure}

\begin{figure}[H]
    \centering
    \includegraphics[width=0.8\linewidth]{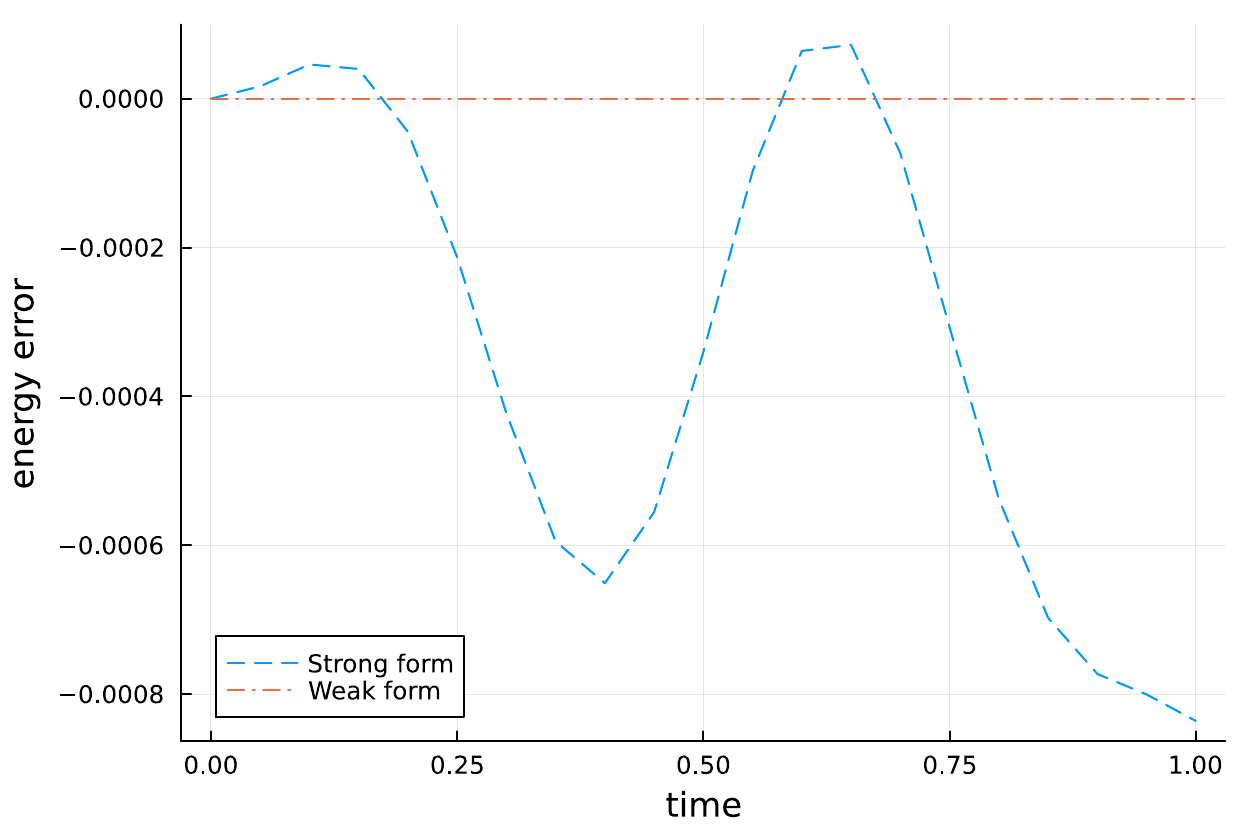}
    \caption{Energy errors for the periodic problem for the non-SBP operator in strong and weak form with SSPRK time integration, $2$ elements with $12$ points per dimension. The maximum absolute energy error for the weak form due to time integration is ~$1.5 \cdot 10^{-13}$.}
    \label{fig:plot_energy_non_sbp}
\end{figure}

\begin{figure}[H]
    \centering
    \includegraphics[width=0.8\linewidth]{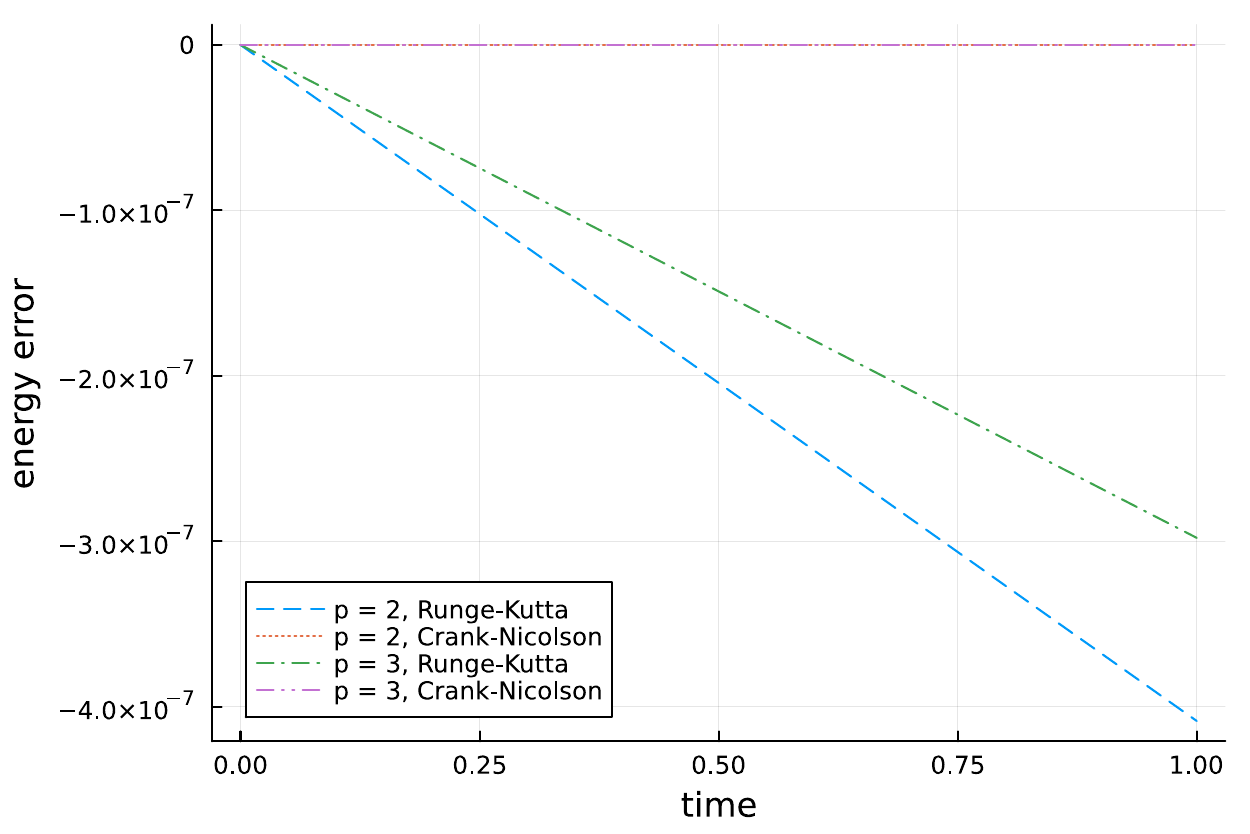}
    \caption{Energy errors for the non-periodic problem for operators of degree $p = 2,3$ with both the SSPRK and Crank-Nicolson time integration, $5$ elements with $20$ points per dimension}
    \label{fig:plot_energy_fine_grid_np}
\end{figure}

We repeat the tests performed for divergence preservation to test the conservation of energy over time.
We use the same operators and the same time integration schemes with the same grid sizes.
In figures \ref{fig:plot_energy_fine_grid} and \ref{fig:plot_energy_fine_grid_np} we plot the difference between the discrete energy at $T = 0$ and the energy at the given point in time.
We can observe that the SSPRK is energy dissipative and does not conserve the energy exactly.
This effect lessens for the finer grid size, but is still present.
On the other hand the Crank Nicolson time integration method conserves the energy exactly as proven previously.
We also run an additional test for the periodic case for energy conservation of the non-SBP operator with mass matrix from the $p = 2$ operator and the derivative matrix from the $p = 3$ operator for both the strong and the weak form of the scheme using $2$ elements and $12$ points per element in each dimension.
To make the influence of the time integration scheme negligible, we use the a time step of $\Delta t = 2\cdot 10^{-5}$ with SSPRK time integration.
Figure \ref{fig:plot_energy_non_sbp} shows that the weak form preserves energy, but the strong form does not do so.
The maximum absolute error for the weak form due to time integration is about $1.5 \cdot 10^{-13}$.
This suggests that for non-SBP operators the weak form is energy conservative, but has a worse convergence behavior, while the strong form has expected convergence rates but is not energy-conservative.

\begin{table}[H]
    \centering
    \begin{tabular}{c|c|c|c|c}
         Number of Elements & $E^x/E^y$ Error & $E^x/E^y$ EOC & $B^z$ Error & $B^z$ EOC   \\
         1 & $1.566639e-2$ & - & $2.98044e-2$ & -\\
         2 & $2.159529e-3$ & $2.86$ & $5.705334e-3$ & $2.39$ \\
         4 & $2.37442e-4$ & $3.19$ & $6.683885e-5$ & $6.42$ \\
         8 & $1.55675e-5$ & $3.93$ & $5.264376e-6$ & $3.67$\\
         16 & $1.795338e-6$ & $3.12$ & $5.574813e-7$ & $3.24$\\
         32 & $2.215563e-7$ & $3.02$ & $2.440927e-8$ & $4.51$ 
    \end{tabular}
    \caption{Discrete $L^2$-Errors and EOC for the periodic problem for degree $p = 3$ with $12$ points per direction per element}
    \label{tbl:eoc_p3_12n}
\end{table}

\begin{table}[H]
    \centering
    \begin{tabular}{c|c|c|c|c}
         Number of Elements & $E^x/E^y$ Error & $E^x/E^y$ EOC & $B^z$ Error & $B^z$ EOC   \\
         1 & $1.487912e-2$ & - & $2.082052e-2$ & -\\
         2 & $2.703255e-3$ & $2.46$ & $5.469707e-4$ & $5.25$ \\
         4 & $1.58823e-4$ & $4.09$ & $4.118751e-5$ & $3.73$ \\
         8 & $7.39334e-6$ & $4.43$ & $5.613232e-6$ & $2.88$\\
         16 & $4.856478e-7$ & $3.93$ & $8.610113e-8$ & $6.03$\\
         32 & $2.558402e-8$ & $4.25$ & $2.16185e-8$ & $1.99$ 
    \end{tabular}
    \caption{Discrete $L^2$-Errors and EOC for the periodic problem for degree $p = 3$ with $13$ points per direction per element}
    \label{tbl:eoc_p3_13n}
\end{table}

\begin{table}[H]
    \centering
    \begin{tabular}{c|c|c|c|c}
         Points per Element & $E^x/E^y$ Error & $E^x/E^y$ EOC & $B^z$ Error & $B^z$ EOC   \\
         12 & $1.566639e-2$ & - & $2.98044e-2$ & -\\
         24 & $1.16444e-3$ & $3.75$ & $3.248053e-3$ & $3.2$ \\
         48 & $9.278422e-5$ & $3.65$ & $1.437682e-4$ & $4.5$ \\
         96 & $5.909218e-6$ & $3.97$ & $6.813822e-6$ & $4.4$\\
         192 & $3.838066e-7$ & $3.94$ & $2.919607e-7$ & $4.54$\\
         384 & $2.347673e-8$ & $4.03$ & $1.309507e-8$ & $4.48$\\
    \end{tabular}
    \caption{Discrete $L^2$-Errors and EOC for the periodic problem for degree $p = 3$ with SSPRK time integration with varying points per element with one element in total}
    \label{tbl:eoc_p3_ssprk}
\end{table}

\begin{table}[H]
    \centering
    \begin{tabular}{c|c|c|c|c}
         Points per Element & $E^x/E^y$ Error & $E^x/E^y$ EOC & $B^z$ Error & $B^z$ EOC   \\
         12 & $1.541206e-2$ & - & $3.321995e-2$ & -\\
         24 & $1.284108e-3$ & $3.59$ & $3.362758e-3$ & $3.3$ \\
         48 & $9.654045e-5$ & $3.73$ & $1.479463e-4$ & $4.51$ \\
         96 & $6.080192e-6$ & $3.99$ & $7.090472e-6$ & $4.38$\\
         192 & $3.978717e-7$ & $3.93$ & $3.070251e-7$ & $4.53$\\
         384 & $2.391862e-8$ & $4.06$ & $1.370273e-8$ & $4.49$\\
    \end{tabular}
    \caption{Discrete $L^2$-Errors and EOC for the periodic problem for the non-SBP operator in the strong form with SSPRK time integration with varying points per dimension per element with one element in total}
    \label{tbl:eoc_non_sbp_strong_ssprk}
\end{table}

\begin{table}[H]
    \centering
    \begin{tabular}{c|c|c|c|c}
         Points per Element & $E^x/E^y$ Error & $E^x/E^y$ EOC & $B^z$ Error & $B^z$ EOC   \\
         12 & $1.840316e-1$ & - & $1.938346e-1$ & -\\
         24 & $9.051131e-2$ & $1.02$ & $9.118747e-2$ & $1.09$ \\
         48 & $2.684118e-2$ & $1.75$ & $1.223276e-2$ & $2.9$ \\
         96 & $9.0215683e-3$ & $1.57$ & $2.224402e-3$ & $2.46$\\
         192 & $3.3507205e-3$ & $1.43$ & $3.815636e-4$ & $2.54$\\
         384 & $1.1464041e-3$ & $1.55$ & $4.608534e-5$ & $3.05$\\
    \end{tabular}
    \caption{Discrete $L^2$-Errors and EOC for the periodic problem for the non-SBP operator in the weak form with SSPRK time integration with varying points per element with one element in total}
    \label{tbl:eoc_non_sbp_weak_ssprk}
\end{table}

\begin{table}[H]
    \centering
    \begin{tabular}{c|c|c|c|c|c|c}
         {$N_e$} & $E^x$ Error & $E^x$ EOC & $E^y$ Error & $E^y$ EOC & $B^z$ Error & $B^z$ EOC   \\
         1 & $1.6079e-3$ & - & $1.8418e-3$ & - & $2.8864e-3$ & -\\
         2 & $1.679e-4$ & $3.26$ & $1.679e-4$ & $3.46$ & $4.7262e-5$ & $5.93$ \\
         4 & $1.1008e-5$ & $3.93$ & $1.1008e-5$ & $3.93$ & $3.7225e-6$ & $3.67$ \\
         8 & $1.2695e-6$ & $3.12$ & $1.2695e-6$ & $3.12$ & $3.942e-7$ & $3.24$\\
         16 & $1.5666e-7$ & $3.02$ & $1.5666e-7$ & $3.02$ & $1.726e-8$ & $4.51$\\
         32 & $1.9483e-8$ & $3.01$ & $1.9483e-8$ & $3.01$ & $7.5131e-10$ & $4.52$ 
    \end{tabular}
    \caption{Discrete $L^2$-Errors and EOC for the non-periodic problem for degree $p = 3$ with $12$ points per direction per element. {Here $N_e$ refers to the number of elements per direction}}
    \label{tbl:eoc_p3_12n_np}
\end{table}

\begin{table}[H]
    \centering
    \begin{tabular}{c|c|c|c|c|c|c}
         Points & $E^x$ Error & $E^x$ EOC & $E^y$ Error & $E^y$ EOC & $B^z$ Error & $B^z$ EOC  \\
         12 & $1.6079e-3$ & - & $1.8418e-3$ & - & $2.8864e-3$ & -\\
         24 & $7.7133e-5$ & $4.38$ & $5.4729e-5$ & $5.07$ & $1.6725e-4$ & $4.11$ \\
         48 & $5.0039e-6$ & $3.95$ & $3.7549e-6$ & $3.87$ & $7.6925e-6$ & $4.44$ \\
         96 & $3.136e-7$ & $4.0$ & $2.8876e-7$ & $3.7$ & $3.5209e-7$ & $4.45$\\
         192 & $1.9237e-8$ & $4.03$ & $2.0145e-8$ & $3.84$ & $1.5232e-8$ & $4.53$\\
         384 & $1.1712e-9$ & $4.04$ & $1.2598-9$ & $4.0$ & $7.5632e-10$ & $4.33$\\
    \end{tabular}
    \caption{Discrete $L^2$-Errors and EOC for the non-periodic problem for degree $p = 3$ with SSPRK time integration with varying points per element per dimension with one element in total}
    \label{tbl:eoc_p3_ssprk_np}
\end{table}

\section{Conclusion}

Using the theory of Finite Element Exterior Calculus \cite{arnold2006,arnold2010}, we derived a family of divergence preserving and semi-discretely energy-conservative schemes for the homogeneous transverse electric Maxwell's equations using pre-existing SBP-FD operators.
Using the Summation-by-Parts property we derived a strong form formulation of our semi-discretization from the weakly discretized Faraday's law, analogously to other SBP schemes such as the DGSEM \cite{kopriva2006}.
To derive the scheme, we extended the construction by Gerritsma et al. \cite{gerritsma2011,preprint_gerritsma2011} using nodal and integral degrees of freedom.
Instead of permitting only polynomial basis functions however, we derived the compatible spaces using more general basis functions, that we only require to have the Lagrange property and the correct derivative values on the grid nodes of the SBP operator. This would allow us to substitute any SBP-FD operator to derive our SBP-FDEC scheme.
Since we extended our approach to 2D and 3D using a tensor product ansatz, we were able to derive general representations for the derivative operators of the 2D and 3D de Rham complexes, yielding commuting discrete complexes, for any given SBP-FD operators, similar to \cite{gerritsma2011}, making the scheme easily extendable to 3D.

Using standard energy-conservative time integration schemes, we can use the semi-discrete energy conservation to derive a fully discrete energy-conservative scheme.
Further we were able to prove the semi-discrete and discrete divergence preserving and the energy conservative properties of our schemes.
We validated our approach with a periodic test case, using SBP-FD operators by Strand \cite{strand1994}, demonstrating both the expected order of convergence, the discrete divergence preserving property and for the Crank-Nicolson time integration method the discrete energy-conservation of our schemes.

\section*{Declarations}

\paragraph{Conflicts of Interest} The authors declare that there are no conflicts of interest.

\paragraph{CrediT authorship contribution statement} 

Daniel Bach: Investigation, Methodology, Formal Analysis, Software, Writing – original draft, Writing – review \& editing. Andrés M. Rueda-Ramírez: Investigation, Methodology, Software, Writing – review \& editing. Eric Sonnendrücker: Writing – review \& editing. David C. Del Rey Fernandez: Writing – review \& editing. Gregor Gassner: Conceptualization, Investigation, Writing – original draft, Writing – review \& editing, Funding acquisition, Project administration.

\paragraph{Data Availability} A reproducibility repository can be accessed under \url{https://github.com/amrueda/paper_2025_sbp_fdec}

\paragraph{Funding}
Gregor J. Gassner and Andrés M. Rueda-Ramírez acknowledge funding through the Klaus-Tschira Stiftung via the project ``HiFiLab'' (00.014.2021) and through the German Federal Ministry for Education and Research (BMBF) project ``ICON-DG'' (01LK2315B) of the  ``WarmWorld Smarter'' program.

Gregor J. Gassner and Daniel Bach acknowledge funding from the German Research Foundation DFG through the research unit ``SNuBIC'' (DFG-FOR5409).

Andrés M. Rueda-Ramírez gratefully acknowledges funding from the Spanish Ministry of Science, Innovation, and Universities through the ``Beatriz Galindo'' grant BG23-00062 and from the European Research Council through the Synergy Grant Agreement No. 101167322-TRANSDIFFUSE.
\appendix

\section{Proof for the vanishing cohomology spaces}

We first note that for $m$ elements in 1D the element boundaries do not change the structure of the global difference matrix from the single element case \eqref{eq:delta_matrix}.
Only the number of rows and columns changes to $m \times N$ and $m \times N + 1$ respectively.
In 2D and 3D we have to account for the fact that the vector of the degrees of freedom cannot be sorted simply from left to right any longer, but this does not affect the component-wise equations if we have a rectangular domain.
Therefore we can restrict ourselves to the single-element case without loss of generality. 
However we have to allow that the element has different amounts of points per direction.
Useful for the proof is the following observation:
\begin{lemma}
    For a vector $\textbf{d} \in \mathbb{R}^N$ and vectors $\textbf{u}^1, \textbf{u}^2 \in \mathbb{R}^{N+1}$ fulfilling
    \begin{equation}
        \label{eq:app_direction_difference}
        u^k_{l} - u^k_{l - 1} = d_l \qquad \forall l \in \{1, \dots N\}, k \in \{1,2\},
    \end{equation}
    it holds that
    \begin{equation}
        \label{eq:app_direction_difference_result}
        u^1_{l} - u^2_{l} = c \qquad \forall l \in \{0, \dots N\},
    \end{equation}    
    for some $c \in \mathbb{R}$.
    Here $\textbf{u}$ has indices from $0$ to $N$ for consistency with the numbering for nodal degrees of freedom.
    Further $u_{l} = \sum_{j=1}^{l} d_i$ is a solution of the inhomogeneous system \eqref{eq:app_direction_difference} and  all solutions have the form $u_{l} = \sum_{j=1}^{l} d_j + c$.
    \\
    
    \textbf{Proof:}
    Equation \eqref{eq:app_direction_difference} is equivalent to
    \begin{equation}
        \label{eq:app_inhomogenous_difference}
        \underline{\Delta} \textbf{u}^k = \textbf{d}, \qquad k \in \{1,2\}.
    \end{equation}
    The matrix $\underline{\Delta}$ has rank $N$, which can easily be seen by computing the determinant of $\underline{\Delta}_{1:N,1:N}$, and its kernel consists of the constant vectors. 
    This proves the first part.
    That $u_{l} = \sum_{j=1}^{l} d_j$ is a solution to \eqref{eq:app_inhomogenous_difference} can be proven by explicit calculation and with the first part every solution has the form $u_l = \sum_{j=1}^{l} d_l + c$.
\end{lemma}

In the following, the number of points in x-direction will be $N_x + 1$, in y-direction $N_y + 1$, and in 3D in z-direction $N_z + 1$.
To prove {$\text{dim}\left(H^0_h\right) = 0$}, we determine the kernel of $\textbf{\underline{curl}}$. We observe that for the first component of $\text{\textbf{\underline{curl}}}(\textbf{u})$ we have for the DOFs after multiplying with $-1$:
\begin{equation}
    \label{eq:app_homogenous_difference}
     u_{i,j} - u_{i,j-1} = 0, \qquad \forall i \in \{0,\dots,N_x\},j \in \{1,\dots,N_y\}.
\end{equation}
Due to lemma 1, all solutions to \eqref{eq:app_homogenous_difference} have the form $u_{i,j} = c_i$ for all $i \in \{0,\dots,N_x\}$.
Inserting this into the equation for the second component we get
\begin{equation}
     c_i - c_{i-1} = 0, \qquad \forall i \in \{1,\dots,N_x\},
\end{equation}
resulting in $u_{i,j} = c$ for some $c \in \mathbb{R}$ for all $i \in \{0,\dots,N_x\}$ and $j \in \{0,\dots,N_y\}$. \\
{This proves the assertion for this cohomology space and arguments for the 2D and 3D gradients are analogous.
For $H^1_h$ }we have to show $\text{Im}\left(\textbf{curl}\right) = \text{Ker}\left(\text{div}\right)$. 
We know that $\text{Im}\left(\underline{\textbf{curl}}\right) \subset \text{Ker}\left(\underline{\text{div}}\right)$ from \eqref{eq:2D_curl_div_identity} and only have to prove the other inclusion.
We also know that $\text{dim}\left(W^0\right) = (N_x+1)(N_y + 1)$ and therefore
\begin{equation}
    \begin{split}
    &\text{dim}\left(\text{Im}\left(\underline{\textbf{curl}}\right)\right) = \text{dim}(W^0) - \text{dim}\left(\text{Ker}\left(\underline{\textbf{curl}}\right)\right) \\
    = \hspace{0.1cm} &(N_x+1)(N_y + 1) - 1 = N_x N_y + N_x + N_y.
    \end{split}
\end{equation} 
We only have to show $\text{dim}\left(\text{Ker}\left(\underline{\text{div}}\right)\right) = N_x N_y + N_x + N_y$.\\
The kernel equation reads
\begin{equation}
    \underline{\Delta}^{2D}_x \textbf{v}^x + \underline{\Delta}^{2D}_y \textbf{v}^y = \textbf{0},
\end{equation}
which is equivalent to
\begin{equation}
    \label{eq:div_component}
    v^x_{i,j} - v^x_{i-1,j} = v^y_{i,j-1} - v^y_{i,j}, \qquad \forall i \in \{1,\dots,N_x\}, j \in \{1,\dots,N_y\} 
\end{equation}
Using lemma 1 with $d_{i,j} \vcentcolon = v^y_{i,j-1} - v^y_{i,j}$, where the index $i$ in \eqref{eq:div_component} corresponds to index $l$ in lemma 1, we obtain $v^x_{i,j} = \sum_{k=1}^{i} \left(v^y_{k,j-1} - v^y_{k,j}\right)  + c_j$ with $c_j \in \mathbb{R}$ for all $j \in \{1,\dots,N_y\}$ and $i \in \{0,\dots,N_x\}$.
Every solution has this form and every choice of $\textbf{v}^y$ and coefficients $c_j$ yields a valid solution.
Therefore 
\begin{equation}
    \text{dim}\left(\text{Ker}\left(\underline{\text{div}}\right)\right) = N_x(N_y+1) + N_y = N_x N_y + N_x + N_y,
\end{equation}
proving { the assertion.
The proof for $\tilde{H}^1_h$ proceeds analogously. \\ 
To prove the claim for $H^2_h$ we have to show $\text{Im}\left(\text{\underline{div}}\right) = W^2$.}
We observe that 
\begin{equation}
    \begin{split}
    &\text{dim}\left(\text{Im}\left(\underline{\text{div}}\right)\right) = \text{dim}(W^1) - \text{dim}\left(\text{Ker}\left(\underline{\text{div}}\right)\right) \\
    = \hspace{0.1cm} &(N_x+1)N_y + N_x(N_y+1) - N_x N_y - N_x - N_y = N_x N_y \\
    = \hspace{0.1cm} &N_x N_y = \text{dim}\left(W^2\right),
    \end{split}
\end{equation}
{ proving the assertion.
Analogously one can prove the same for $\tilde{H}^2_h$.\\
This leaves us with the discrete cohomology spaces $H^{1}_{3D,h}$, $H^{2}_{3D,h}$, and $H^{3}_{3D,h}$.
To prove $\text{dim}\left(H^{1}_{3D,h}\right) = 0$} we need to prove $\text{Im}\left(\underline{\textbf{grad}}\right) = \text{Ker}\left(\underline{\textbf{curl}}\right)$.
Since $\text{dim}\left(\text{Ker}\left(\underline{\textbf{grad}}\right)\right) = 1$, we have 
\begin{equation}
    \begin{split}&\text{dim}\left(\text{Im}\left(\underline{\textbf{grad}}\right)\right) = \text{dim}\left(W^0\right) - \text{dim}\left(\text{Ker}\left(\underline{\textbf{grad}}\right)\right) \\
    = \hspace{0.1cm}& (N_x+1)(N_y + 1)(N_z + 1) - 1 \\
    = \hspace{0.1cm}& N_x N_y N_z + N_x N_y + N_y N_z + N_z N_x + N_x + N_y + N_z.
    \end{split}
\end{equation}
The equations for the $\text{Ker}\left(\underline{\textbf{curl}}\right)$ are
\begin{align}
    \label{eq:app_curl_kernel_1}
    &\underline{\Delta}^{3D}_y \textbf{v}^z - \underline{\Delta}^{3D}_z \textbf{v}^y = \textbf{0}, \\
    \label{eq:app_curl_kernel_2}
    &\underline{\Delta}^{3D}_z \textbf{v}^x - \underline{\Delta}^{3D}_x \textbf{v}^z = \textbf{0}, \\
    \label{eq:app_curl_kernel_3}
    &\underline{\Delta}^{3D}_x \textbf{v}^y - \underline{\Delta}^{3D}_y \textbf{v}^x = \textbf{0},
\end{align}
where $\underline{\Delta}^{3D}$ matrices are defined similar to the 2D ones, but with a triple Kronecker product with two of the arguments being the identity matrix. 
Similar to the argument for { $H^1_h$} we obtain from equations \eqref{eq:app_curl_kernel_1} and \eqref{eq:app_curl_kernel_2}
\begin{align}
    v^x_{i,j,k} &= \sum_{l=1}^{k} \left(v^z_{i,j,l} - v^z_{i-1,j,l}\right) + c_{i,j},\\
    v^y_{i,j,k} &= \sum_{l=1}^{k} \left(v^z_{i,j,l} - v^z_{i,j-1,l}\right) + d_{i,j},
\end{align}
where $v^x_{i,j,k}$ has indices $i \in \{1, \dots, N_x\}$, $j \in \{0, \dots, N_y\}$, $k \in \{0, \dots, N_z\}$, and $v^y_{i,j,k}$ has indices $i \in \{0, \dots, N_x\}$, $j \in \{1, \dots, N_y\}$, and $k \in \{0, \dots, N_z\}$.
Inserting this in \eqref{eq:app_curl_kernel_3}, we get component-wise
\begin{equation}
    \begin{split}
    &\sum_{l=1}^{k} \left(v^z_{i,j,l} - v^z_{i,j-1,l} - v^z_{i-1,j,l} + v^z_{i-1,j-1,l}\right) + \left(c_{i,j} - c_{i,j-1}\right) \\
    = &\sum_{l=1}^{k} \left(v^z_{i,j,l} - v^z_{i,j-1,l} - v^z_{i-1,j,l} + v^z_{i-1,j-1,l}\right) + \left(d_{i,j} - d_{i-1,j}\right)
    \end{split}
\end{equation}
for all $i\in \{1,\dots,N_x\}$, $j\in \{1,\dots,N_y\}$, and $k \in \{0,\dots,N_z\}$.
This is equivalent to
\begin{equation}
    \left(c_{i,j} - c_{i,j-1}\right) = \left(d_{i,j} - d_{i-1,j}\right), \qquad \forall i\in \{1,\dots,N_x\}, j\in \{1,\dots,N_y\},
\end{equation}
and thus we obtain
\begin{equation}
    c_{i,j} = \sum_{l = 1}^j \left(d_{i,l} - d_{i-1,l}\right) + e_i, \qquad \qquad \forall i \in \{1,\dots,N_x\}, j \in \{0,\dots,N_y\}.
\end{equation}
Counting degrees of freedom we have $(N_x+1) (N_y + 1) N_z$ from $\textbf{v}^z$, $(N_x+1)N_y$ from the coefficients $d_{i,j}$ and $N_x$ from the coefficients $e_i$.
Together we have 
\begin{equation}
    \begin{split}
    &(N_x+1) (N_y + 1) N_z + (N_x+1)N_y + N_x \\
    = \hspace{0.1cm} &N_x N_y N_z + N_x N_y + N_y N_z + N_x N_z + N_x + N_y + N_z \\
    = \hspace{0.1cm} &(N_x+1)(N_y + 1)(N_z + 1) - 1 = \text{dim}\left(\text{Im}\left(\underline{\textbf{grad}}\right)\right),
    \end{split}
\end{equation}
which proves {$\text{dim}\left(H^{1}_{3D,h}\right) = 0$}.\\
Since $\text{dim}\left(W^1\right) = N_x (N_y + 1) (N_z + 1) + (N_x + 1) N_y (N_z + 1) + (N_x + 1) (N_y + 1) N_z$, we have 
\begin{equation}
    \begin{split}
    &\text{dim}\left(\text{Im}\left(\underline{\textbf{curl}}\right)\right) = \text{dim}\left(W^1\right) - \text{dim}\left(\text{Ker}\left(\underline{\textbf{curl}}\right)\right) \\
    = \hspace{0.1cm}&N_x (N_y + 1) (N_z + 1) + (N_x + 1) N_y (N_z + 1) + (N_x + 1) (N_y + 1) N_z\\
    &- N_x N_y N_z - N_x N_y - N_y N_z - N_x N_z - N_x - N_y - N_z \\
    = \hspace{0.1cm}&2N_x N_y N_z + N_x N_y + N_y N_z + N_z N_x.
    \end{split}
\end{equation}
{For $H^2_{3D,h}$ we have the kernel equation for the \underline{div} operator}
\begin{equation}
    \underline{\Delta}^{3D}_x \textbf{w}^x + \underline{\Delta}^{3D}_y \textbf{w}^y +\underline{\Delta}^{3D}_z \textbf{w}^z = \textbf{0},
\end{equation}
resulting in 
\begin{equation}
    w^x_{i,j,k} = \sum_{l = 1}^i \left(\textbf{w}^y_{l,j-1,k} - \textbf{w}^y_{l,j,k} + \textbf{w}^z_{l,j,k-1} - \textbf{w}^z_{l,j,k}\right) + f_{j,k},
\end{equation}
for all $i \in \{0, \dots, N_x\}$, $j \in \{1, \dots, N_y\}$, and $k \in \{1, \dots, N_z\}$.
Counting degrees of freedom again we get $N_x (N_y + 1) N_z$ from $\textbf{w}^y$, $N_x N_y (N_z + 1)$ from $\textbf{w}^z$, and $N_y N_z$ from the coefficients $f_{j,k}$, resulting in 
\begin{equation}
    \begin{split}
    &N_x (N_y + 1) N_z + N_x N_y (N_z + 1) + N_y N_z \\
    = \hspace{0.1cm} &2N_x N_y N_z + N_x N_y + N_y N_z + N_z N_x = \text{dim}\left(\text{Im}\left(\underline{\textbf{curl}}\right)\right),
    \end{split}
\end{equation}
proving {$\text{dim}\left(H^2_{3D,h}\right) = 0$.
For the last equation, we note that}
\begin{equation}
    \begin{split}&
        \text{dim}\left(\text{Im}\left(\underline{\text{div}}\right)\right) = \text{dim}\left(W^2\right) - \text{dim}\left(\text{Ker}\left(\underline{\text{div}}\right)\right) \\
        = \hspace{0.1cm}& (N_x + 1) N_y N_z + N_x (N_y + 1) N_z + N_x N_y (N_z + 1) \\
        &- 2N_x N_y N_z - N_x N_y - N_y N_z - N_z N_x \\
        = \hspace{0.1cm}& N_x N_y N_z = \text{dim}\left(W^3\right),
    \end{split}
\end{equation} which completes the proof.

\section{Further test results}

\begin{table}[H]
    \centering
    \begin{tabular}{c|c|c|c|c}
         Number of Elements & $E^x/E^y$ Error & $E^x/E^y$ EOC & $B^z$ Error & $B^z$ EOC   \\
         1 & $1.487802e-1$ & - & $2.128334e-1$ & -\\
         2 & $1.457635e-2$ & $3.35$ & $3.770852e-2$ & $2.50$ \\
         4 & $3.215945e-3$ & $2.18$ & $2.840422e-3$ & $3.73$ \\
         8 & $4.087012e-4$ & $2.98$ & $2.402008e-4$ & $3.56$\\
         16 & $4.357807e-5$ & $3.23$ & $4.276273e-5$ & $2.49$\\
         32 & $5.590312e-6$ & $2.96$ & $8.925102e-6$ & $2.26$ 
    \end{tabular}
    \caption{Discrete $L^2$-Errors and EOC for the periodic problem for degree $p = 2$ with $8$ points per direction per element}
    \label{tbl:eoc_p2_8n}
\end{table}

\begin{table}[H]
    \centering
    \begin{tabular}{c|c|c|c|c}
         Number of Elements & $E^x/E^y$ Error & $E^x/E^y$ EOC & $B^z$ Error & $B^z$ EOC   \\
         1 & $9.964855e-2$ & - & $1.260615e-1$ & -\\
         2 & $1.130372e-2$ & $3.14$ & $5.891993e-3$ & $4.42$ \\
         4 & $2.212383e-3$ & $2.35$ & $1.149555e-3$ & $2.36$ \\
         8 & $5.62475e-4$ & $1.98$ & $1.331211e-4$ & $3.11$\\
         16 & $1.407649e-4$ & $2$ & $5.033169e-6$ & $4.73$\\
         32 & $3.518716e-5$ & $2$ & $1.586997e-6$ & $1.67$ 
    \end{tabular}
        \caption{Discrete $L^2$-Errors and EOC for the periodic problem for degree $p = 2$ with $9$ points per direction per element}
    \label{tbl:eoc_p2_9n}
\end{table}

\begin{table}[!h]
    \centering
    \begin{tabular}{c|c|c|c|c}
         Number of Elements & $E^x/E^y$ Error & $E^x/E^y$ EOC & $B^z$ Error & $B^z$ EOC   \\
         1 & $4.418904e-2$ & - & $1.054903e-1$ & -\\
         2 & $3.740646e-3$ & $3.56$ & $1.440109e-2$ & $2.87$ \\
         4 & $1.34654e-3$ & $1.47$ & $1.170316e-3$ & $3.62$ \\
         8 & $1.722958e-4$ & $2.97$ & $2.350005e-4$ & $2.32$\\
         16 & $2.678038e-5$ & $2.69$ & $1.035082e-5$ & $4.5$\\
         32 & $2.699576e-6$ & $3.31$ & $2.375989e-6$ & $2.12$ 
    \end{tabular}
    \caption{Discrete $L^2$-Errors and EOC for the periodic problem for degree $p = 2$ with $10$ points per direction per element}
    \label{tbl:eoc_p2_10n}
\end{table}

\begin{table}[!h]
    \centering
    \begin{tabular}{c|c|c|c|c}
         Number of Elements & $E^x/E^y$ Error & $E^x/E^y$ EOC & $B^z$ Error & $B^z$ EOC   \\
         1 & $8.506095e-3$ & - & $2.235792e-2$ & -\\
         2 & $9.643637e-4$ & $3.14$ & $3.285374e-3$ & $2.77$ \\
         4 & $1.017081e-4$ & $3.25$ & $4.232642e-5$ & $6.28$ \\
         8 & $7.697999e-6$ & $3.72$ & $4.7801012e-6$ & $3.15$\\
         16 & $9.530514e-7$ & $3.01$ & $1.958355e-7$ & $4.61$\\
         32 & $1.176027e-7$ & $3.02$ & $1.007924e-8$ & $4.28$ 
    \end{tabular}
    \caption{Discrete $L^2$-Errors and EOC for the periodic problem for degree $p = 3$ with $14$ points per direction per element}
    \label{tbl:eoc_p3_14n}
\end{table}

\begin{table}[H]
    \centering
    \begin{tabular}{c|c|c|c|c}
         Points per Element & $E^x/E^y$ Error & $E^x/E^y$ EOC & $B^z$ Error & $B^z$ EOC   \\
         8 & $1.487802e-1$ & - & $2.128334e-1$ & -\\
         16 & $9.078998e-3$ & $4.03$ & $1.586406e-2$ & $3.75$ \\
         32 & $4.807838e-4$ & $4.24$ & $8.52153e-4$ & $4.22$ \\
         64 & $3.536512e-5$ & $3.76$ & $4.897816e-5$ & $4.12$\\
         128 & $2.885322e-6$ & $3.62$ & $2.921983e-6$ & $4.07$\\
         256 & $2.459519e-7$ & $3.55$ & $1.785489e-7$ & $4.03$\\
    \end{tabular}
    \caption{Discrete $L^2$-Errors and EOC for the periodic problem for degree $p = 2$ with SSPRK time integration with varying points per element per direction with one element in total}
    \label{tbl:eoc_p2_ssprk}
\end{table}

\begin{table}[H]
    \centering
    \begin{tabular}{c|c|c|c|c}
         Points per Element & $E^x/E^y$ Error & $E^x/E^y$ EOC & $B^z$ Error & $B^z$ EOC   \\
         8 & $1.487802e-1$ & - & $2.128334e-1$ & -\\
         16 & $9.078998e-3$ & $4.03$ & $1.586407e-2$ & $3.75$ \\
         32 & $4.807839e-4$ & $4.24$ & $8.521539e-4$ & $4.22$ \\
         64 & $3.536516e-5$ & $3.76$ & $4.897902e-5$ & $4.12$\\
         128 & $2.885341e-6$ & $3.62$ & $2.922832e-6$ & $4.07$\\
         256 & $2.459736e-7$ & $3.55$ & $1.794117e-7$ & $4.03$
    \end{tabular}
    \caption{Discrete $L^2$-Errors and EOC for the periodic problem for degree $p = 2$ with Crank-Nicolson time integration with varying points per element per direction with one element in total}
    \label{tbl:eoc_p2_avf}  
\end{table}

\begin{table}[H]
    \centering
    \begin{tabular}{c|c|c|c|c}
         Points per Element & $E^x/E^y$ Error & $E^x/E^y$ EOC & $B^z$ Error & $B^z$ EOC   \\
         12 & $1.566639e-2$ & - & $2.98044e-2$ & -\\
         24 & $1.16444e-3$ & $3.75$ & $3.248053e-3$ & $3.2$ \\
         48 & $9.278422e-5$ & $3.65$ & $1.437682e-4$ & $4.5$ \\
         96 & $5.909218e-6$ & $3.97$ & $6.813825e-6$ & $4.4$\\
         192 & $3.838071e-7$ & $3.94$ & $2.919745e-7$ & $4.54$\\
         384 & $2.348321e-8$ & $4.03$ & $1.339557e-8$ & $4.45$
    \end{tabular}
    \caption{Discrete $L^2$-Errors and EOC for the periodic problem for degree $p = 3$ with Crank-Nicolson time integration with varying points per element per direction with one element in total}
    \label{tbl:eoc_p3_avf} 
\end{table}

\begin{figure}[H]
    \centering
    \includegraphics[width=0.8\linewidth]{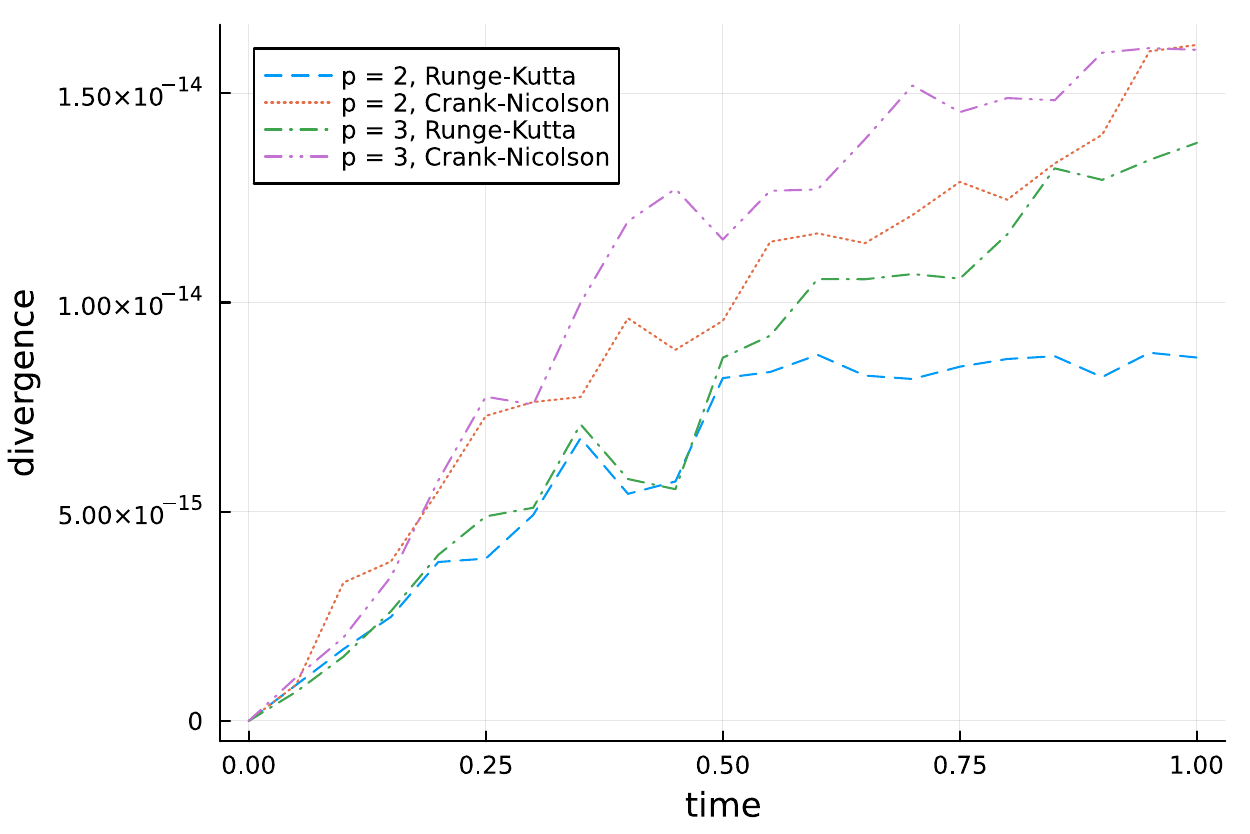}
    \caption{Divergence errors for the periodic problem for operators of degree $p = 2,3$ with both the SSPRK and Crank-Nicolson time integration, $2$ elements with $12$ points per dimension}
    \label{fig:plot_div_coarse_grid}
\end{figure}

\begin{figure}[H]
    \centering
    \includegraphics[width=0.8\linewidth]{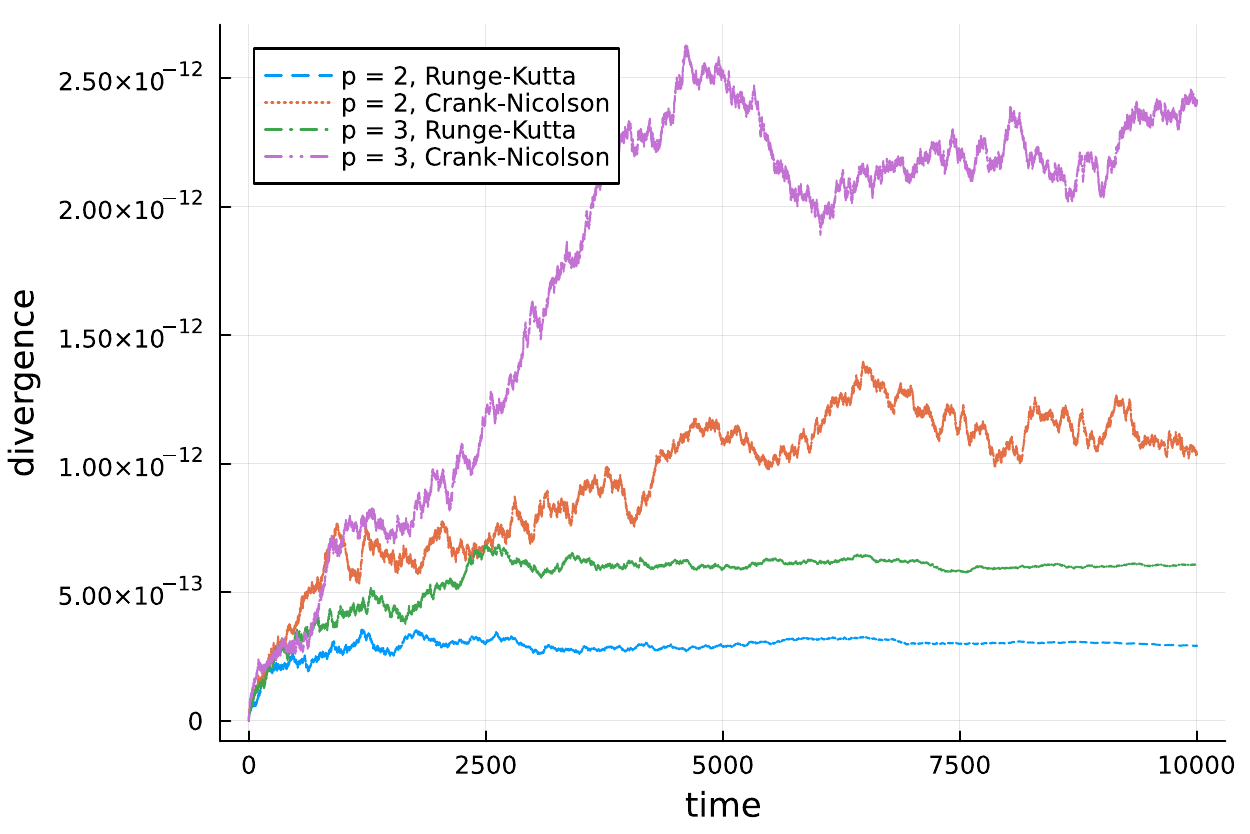}
    \caption{Divergence errors for the periodic problem for operators of degree $p = 2,3$ with both the SSPRK and Crank-Nicolson time integration, $2$ elements with $12$ points per dimension, end time $T = 10000$}
    \label{fig:plot_div_coarse_grid_T10000}
\end{figure}

\begin{figure}[H]
    \centering
    \includegraphics[width=0.8\linewidth]{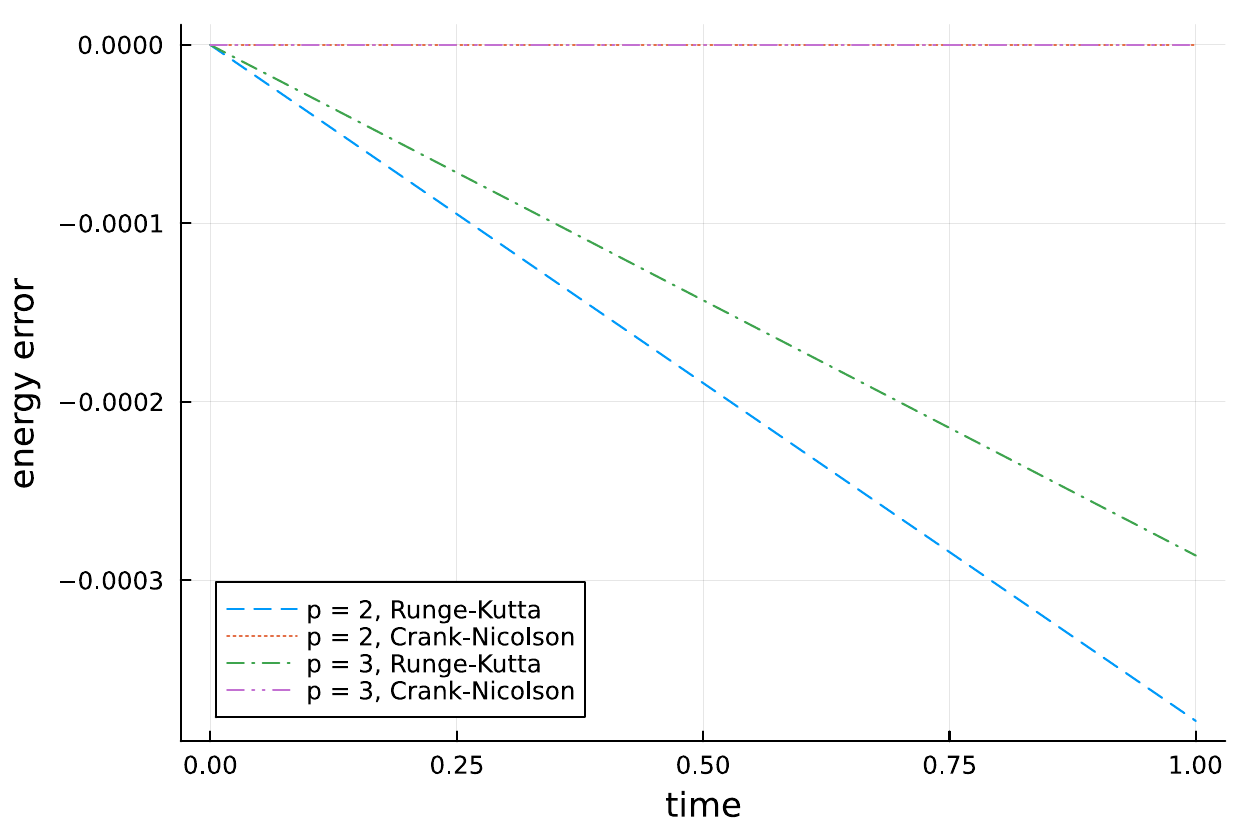}
    \caption{Energy errors for the periodic problem for operators of degree $p = 2,3$ with both the SSPRK and Crank-Nicolson time integration, $2$ elements with $12$ points per dimension}
    \label{fig:plot_energy_coarse_grid}
\end{figure}

\begin{figure}[H]
    \centering
    \includegraphics[width=0.8\linewidth]{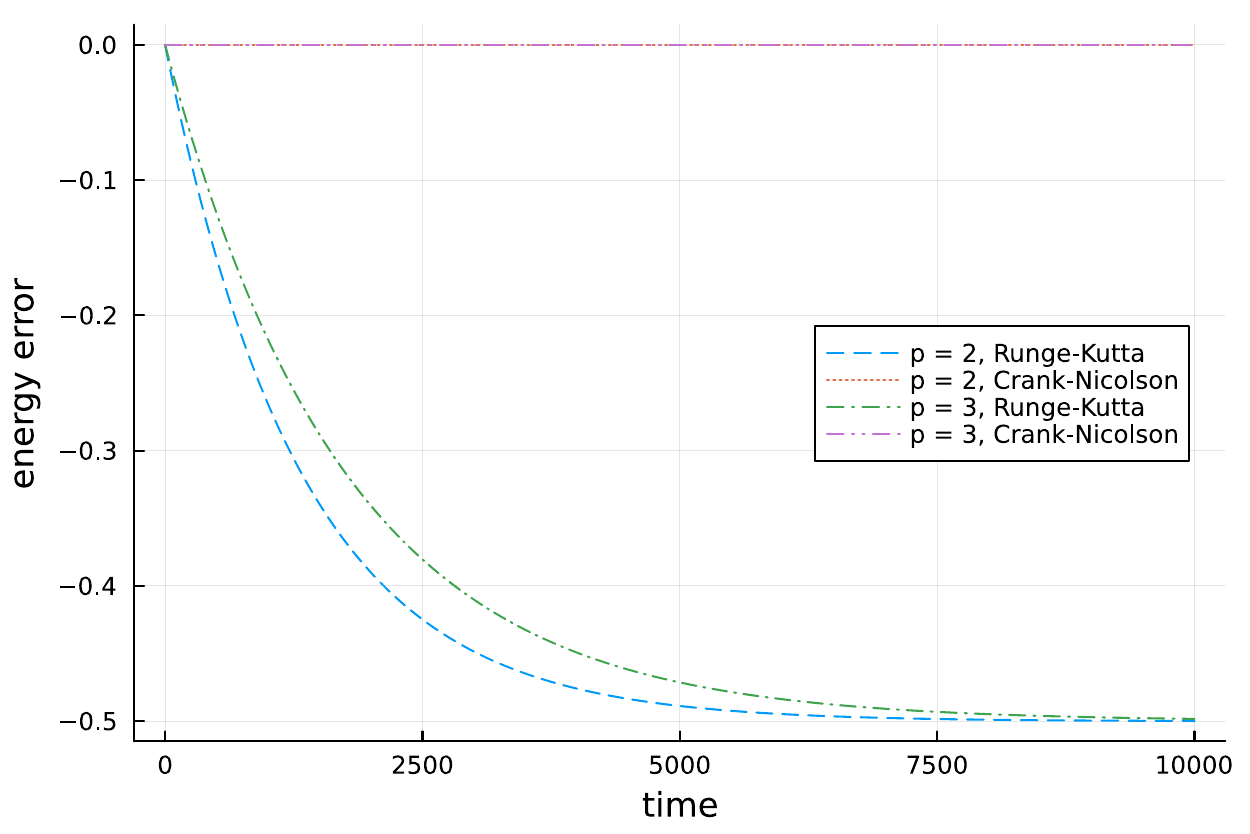}
    \caption{Energy errors for the periodic problem for operators of degree $p = 2,3$ with both the SSPRK and Crank-Nicolson integration, $2$ elements with $12$ points per dimension, end time $T = 10000$}
    \label{fig:plot_energy_coarse_grid_T10000}
\end{figure}

\section{Assembly of the Global Spaces}

To define the global spaces, we assume that the existence of a Lagrange basis $\{\tilde{l}_i^k\}_{i = 0}^N$ as specified in \ref{subsection:basis_functions} for each element $K^k$, that we can get as a linear transformation of such a basis from a reference element.
From these we can derive the histopolation functions $\{\tilde{h}_i^k\}_{i = 1}^N$ in the same way as in the one-element-case.
We construct the global basis functions $\{l_i^k\}_{i = 0, k = 1}^{N-1,m}$ and $\{h_i^k\}_{i,k = 1}^{N,m}$ {analogous to continuous Galerkin finite elements} as follows for all $k \in \{1, \dots, m\}$:
\begin{align}
    \label{eq:global_lagrange_basis_inner}
    l_i^k(x) &:= \begin{cases}
        \tilde{l}_i^k(x) \text{ for }x \in K^k, \\
        0 \text{ for } x \in K^l, l \neq k,
    \end{cases}
    \qquad \forall i \in \{1,\dots,N-1\}, \\
    \label{eq:global_lagrange_basis_boundary_0}
    l_0^k(x) &:= \begin{cases}
        \tilde{l}_0^k(x) \text{ for }x \in K^k, \\
        \tilde{l}_N^{\tilde{\sigma}(k)}(x) \text{ for }x \in {K^{\tilde{\sigma}(k)}} ,\\   
        0, \text{ for } x \in K^l, l \neq k, l \neq {{\tilde{\sigma}(k)}}
    \end{cases} \\
    \label{eq:global_lagrange_basis_boundary_N}
    l_N^k(x) &:= \begin{cases}
        \tilde{l}_N^k(x) \text{ for }x \in K^k, \\
        \tilde{l}_0^{\sigma(k)}(x) \text{ for }x \in K^{\sigma(k)}, \\ 
        0 \text{ for } x \in K^l, l \neq k, l \neq {\sigma(k)},
    \end{cases} \\
    \label{eq_global_histopolation_basis}
    h_i^k(x) &:= \begin{cases}
        \tilde{h}_i^k(x) \text{ for }x \in K^k, \\
        0 \text{ for } x \in K^l, l \neq k,
    \end{cases}
    \qquad \forall i \in \{1,\dots,N\},
\end{align}
where $\sigma(k)$ is defined as the element index of the right neighbor of the element with index $k$ and $\tilde{\sigma}(k)$ is the index of the left neighbor of element $k$.
Here $l_N^k$ and $l_0^{\sigma(k)}$ refer to the same function, as do $l_0^k$ and $l_N^{\tilde{\sigma}(k)}$.
In the non-periodic case such neighbouring elements will not always exist, in which case we set $K^{\sigma(k)}$ or $K^{\tilde{\sigma}(k)}$ respectively to be the empty set.
While the Lagrange basis is globally continuous, the histopolation basis does not have a unique value at element boundaries.
When we use the nodal evaluation matrix for the histopolation functions \eqref{eq:vandermonde}, we will do so for each element separately and obtain the value of the continuous extension of the basis functions to the boundary when restricted to the element in question.
Analogously to \ref{subsection:basis_functions} we can now define the space $U_l$ as the span of the global Lagrange functions and $U_h$ as the span of the global histopolation functions.
We again have $\bigpartialderiv{}{x}(U_l) = U_h$ for the non-periodic case, where $\bigpartialderiv{}{x}$ is the derivative operator.
This is because there is an additional continuity requirement at the domain boundary that breaks down when taking the derivative.In the periodic case we have $\bigpartialderiv{}{x}(U_l) \subset U_h$ instead, since the kernel and the number of DOFs for $U_h$ remains the same, but we lose one DOF for $U_l$.
Now, in analogy to \eqref{eq:def_interpolation} and \eqref{eq:def_histopolation}, let $p^0$ be the global interpolation operator and $p^1$ the global histopolation operator associated with the defined bases.
We can apply the derivative formula \eqref{eq:derivative_formula} element-wise to obtain the commutation property for the multi-element case.
This means the diagram in Fig. \ref{fig:1d_discrete_de_rham} for the one-dimensional de Rham complex commutes.
\bibliography{references}
\bibliographystyle{elsarticle-num}

\end{document}